\documentclass[a4paper,twoside,draft,11pt,righttag]{amsart}
\usepackage{amsxtra}
\usepackage{amsfonts}
\usepackage{amssymb}
\renewcommand{\arraystretch}{.35}

\theoremstyle{plain}
\newtheorem{lema}{Lemma}[section]
\newtheorem{prop}[lema]{Proposition}%[section]
\newtheorem{teo}[lema]{Theorem}%[section]
\newtheorem{coro}[lema]{Corollary}%[section]

\theoremstyle{definition}
\newtheorem{ejem}[lema]{Example}%[section]
\newtheorem{rem}[lema]{Remark}%[section]

\theoremstyle{remark}

\numberwithin{equation}{section}

\newcommand{\Z}{\mathbb Z}

\newcommand{\R}{\mathbb R}
\newcommand{\C}{\mathbb C}

\newcommand{\g}{\gamma}
\newcommand{\ld}{\lambda}
\newcommand{\Ld}{\Lambda}
\newcommand{\G}{\Gamma}

\newcommand{\noi}{\noindent}
\newcommand{\f}{\frac}
\newcommand{\tf}{\tfrac}

\newcommand{\arr}{\rightarrow}

\newcommand{\sk}{\smallskip}
\newcommand{\msk}{\medskip}

\newcommand{\vs}{\vspace{1em}}

\newcommand{\vcp}{\Gamma \backslash \R^n}
\newcommand{\I}{\text{\sl Id}}

\title[Spectral properties of compact flat 4-manifolds]{Spectral Properties of 4-dimensional \\ compact flat manifolds}
\author[R. J. Miatello]{Roberto J. Miatello}
\address{FaMAF--CIEM \\ Universidad Nacional de C\'ordoba \\ Argentina.}
\email{miatello@mate.uncor.edu, podesta@mate.uncor.edu}
\author[R. A. Podest\'a]{Ricardo A. Podest\'a}
\thanks{Partially supported by Conicet and grants from Secyt-UNC, Foncyt and AgCba.}

\begin{document}
\bibliographystyle{plain}

\begin{abstract}
    We study the spectral properties of a large class of compact flat Riemannian manifolds of dimension 4, namely, those  whose corresponding Bieberbach groups have the canonical lattice as translation lattice.
By using the explicit expression of the heat trace of the
Laplacian acting on $p$-forms, we determine all $p$-isospectral
and $L$-isospectral pairs and we show that in this class of
manifolds, isospectrality on functions and isospectrality on
$p$-forms  for all values of $p$ are equivalent to each other. The
list shows  for any $p$, $1\le p \le 3$, many $p$-isospectral
pairs that are not isospectral on functions and have different
lengths of closed geodesics. We also determine all length
isospectral pairs (i.e.\@ with the same length multiplicities),
showing that there are two weak length isospectral pairs that are
not length isospectral, and many pairs, $p$-isospectral for all
$p$ and not length isospectral.
\end{abstract}

\maketitle

\section*{Introduction}                     \label{s.intro}
    If $M$ is a compact Riemannian $n$-manifold, let  $\text{Spec}_p(M)$ denote the $p$-spectrum of $M$, that is,  the collection of eigenvalues, counted with multiplicities,  of the Laplacian acting on smooth $p$-forms. If $\text{Spec}_p(M)=\text{Spec}_p(M')$, then
$M$ and $M'$ are said to be {\em $p$-isospectral}.
In the function case, i.e.\@ if $p=0$, one says that $M$ and $M'$ are {\em isospectral}. The set of all lengths of closed geodesics of $M$ (resp.\@ lengths with multiplicities) is called the {\em weak length spectrum} or {\em $L$-spectrum} (resp.\@ {\em length spectrum} or {\em $[L]$-spectrum}) of $M$.

In \cite{MR2,MR3} multiplicity formulas for the eigenvalues of the
Laplacian acting on natural vector bundles over compact flat
manifolds were given, together with isospectrality criteria, that
allow to give a variety of examples and counterexamples. In the
multiplicity formulas for eigenvalues for the Laplacian acting on
$p$-forms, certain traces appear. For flat manifolds of diagonal
type (see Section 1) these are given by integral values of
Krawtchouk polynomials, hence one can use their integral zeros to
construct examples of $p$-isospectral manifolds that are not
0-isospectral. In \cite{MR4}  different types of length spectra
are studied and  $p$-isospectrality is compared with length
isospectrality. Some of the examples in \cite{MR2, MR3, MR4} occur
already in dimension 4 and, given that there is a full
classification of the corresponding Bieberbach groups in this case
(\cite{BBNWZ}),  it seems natural to study the spectral properties
of flat manifolds  for $n=4$. This is the main goal of the present
paper.

It is known that if two flat manifolds are isospectral, the
corresponding covering tori must be isospectral to each other
(\cite{Su}). It is thus natural, when comparing spectra, to
consider  flat manifolds both having the same covering torus. In
this paper we will restrict ourselves to manifolds covered by the
standard torus, that is, the corresponding Bieberbach groups have
the canonical (cubic) lattice of translations. We will be left
with  54 out of the 74 diffeomorphism classes of four dimensional
manifolds and we shall see that some of the diffeomorphism classes
allow  several different isometry classes. If we include
representatives for all these isometry classes, then we end with a
total of 79 different isometry classes of Bieberbach 4-manifolds,
covered by the standard 4-dimensional torus.

 An outline of the paper is as follows. In Section 1 we briefly review  the main facts
 on Bieberbach groups  and some results from \cite{MR2,MR3,MR4}.  We also give a table with
 the values of the Krawtchouk polynomials $K_p^4(j)$  for $0\leq p,j \leq 4$
 (see (\ref{krawtvalues})).

In Section~\ref{s.sec2} we give a full set of representatives for the
isometry classes in each diffeomorphism class, and compute the ingredients
for the $p$-heat traces and the lengths of closed geodesics for all
Bieberbach groups in our class. At the end of the section, we use formula
(\ref{multip}) to give explicit computations of multiplicities of
eigenvalues in the particular case of two flat manifolds with holonomy
group $D_4$. These turn out to be $p$-isospectral if and only if $p=1,3$.
Furthermore they do not have the same lengths of closed geodesics. This is
the first example with these properties for flat manifolds having {\it
nonabelian} holonomy groups.

In Section \ref{s.sec3} we give expressions for the $p$-heat
traces $Z_p^\G(s)$ in the context of this paper and then rewrite
$Z_p^\G(s)$ in terms of polynomial expressions that are simpler
than the zeta functions, and still encode the same  spectral
information.

In the last section we give all pairs of Sunada-isospectral,
$p$-isospectral and $L$-isospectral flat manifolds of dimension 4
in our context, by comparison of the zeta functions $Z_p^\G(s)$.
We will show that all  0-isospectral pairs are actually
$L$-isospectral, $p$-isospectral for all $p$ and have the same
holonomy representation. Moreover, the non-diffeomorphic pairs
occur only in the class of Bieberbach manifolds of diagonal type.
On the other hand, our list will reveal a large number of
$p$-isospectral pairs for either $p=1$ and $p=3$, or for $p=2$,
that are not isospectral on functions and have different lengths
of closed geodesics. Also, we shall see that there exist several
$L$-isospectral pairs that are not isospectral. Finally, we take
into account length multiplicities, finding all [$L$]-isospectral
pairs, i.e.\@ having the same lengths of closed geodesics and the
same multiplicities for each length. The list obtained shows,
already in dimension 4, that {\em non diffeomorphic} Bieberbach
manifolds that are $p$-isospectral for all $0\le p \le 4$, are
often not [$L$]-isospectral (4 pairs out of 9).

\section{Preliminaries}                     \label{s.prelim}
    We shall first recall some standard facts on flat Riemannian manifolds (see \cite{Ch}).
A discrete, cocompact subgroup  $\Gamma$ of  the isometry group of
$\R^n$,  $I(\R^n)$, is called a crystallographic group. If
furthermore, $\G$ is torsion-free, then $\G$ is said to be a {\it
Bieberbach group}. Such $\Gamma$ acts properly discontinuously on
$\R^n$, thus $M_\Gamma = \Gamma\backslash\R^n$ is a compact flat
Riemannian manifold with fundamental group $\Gamma$. Any such
manifold arises in this way.  Any element $\gamma \in I(\R^n)$
decomposes uniquely as $\gamma = B L_b$, with $B \in \text{O}(n)$,
$b\in \R^n$ and $L_b$ denotes translation by $b$. The translations
in $\Gamma$ form a normal maximal abelian subgroup of finite
index, identified with $\Lambda$, a lattice in $\R^n$ which is
$B$-stable for each $BL_b \in \Gamma$. The quotient $F :=
\Lambda\backslash\Gamma$ is called the holonomy group of $\Gamma$
and gives the linear holonomy group of the Riemannian manifold
$M_\Gamma$. The action of $F$ on $\Lambda$ defines an integral
representation of $F$, usually called the holonomy representation.

 If $BL_b$ is in $\G$, denote $n_B:=\text{dim} \ker(B-\I)$.  The torsion-free condition on $\G$ implies that for any $BL_b \in \G$ we have $n_B>0$ and furthermore, $b_+ :=p_B(b)\ne 0$,
where $p_B$ denotes the orthogonal projection onto $\text{ker}(B-\I)$.

We now recall from \cite{MR2, MR3}  some facts on the spectrum of
Laplace operators on vector bundles over flat manifolds. If $\tau$
is an irreducible representation of $K=\text{O}(n)$ and
$G=I(\R^n)$ we form the vector bundle $E_\tau$ over $G/K \simeq
\R^n$ associated to $\tau$ and consider the corresponding bundle
$\Gamma\backslash E_\tau$ over $\Gamma\backslash  \R^n =M_\Gamma$.
Let $-\Delta_\tau$ be the connection Laplacian on this bundle. For
any  nonnegative real number $\mu$, let $\Lambda^*_\mu=\{\lambda
\in \Lambda^* : \| \lambda\|^2=\mu\}$, where $\Ld^*$ denotes the
dual lattice of $\Ld$.  In \cite{MR3}, Theorem 2.1,
 it is shown that the multiplicity of the eigenvalue $4\pi^2 \mu$ of
$-\Delta_\tau$ is given by
\begin{equation}                                \label{multip}
 d_{\tau,\mu}(\Gamma)= \tfrac 1{|F|} \sum_{\gamma=BL_b \in \Lambda\backslash \Gamma}\text{tr}\,\tau(B) \,e_{\mu,\gamma}
\end{equation}
where $e_{\mu,\gamma} = \sum_{v\in {\Lambda^*_\mu}:Bv=v} e^{-2\pi
i v\cdot b}$. In the case when $\tau=\tau_p$, the $p$-exterior
representation of $\text{O}(n)$, we shall write $\text{tr}_p(B)$
and $d_{p,\mu}(\Gamma)$ in place of $\text{tr}\,\tau_p(B)$ and $
d_{\tau_p,\mu}(\Gamma)$ respectively. If $p=0$ we have
$\Delta_{\tau_0}=\Delta$, the standard Laplacian on functions.

For a special class of flat manifolds the terms in formula (\ref{multip}) can be made more explicit.
A Bieberbach group $\Gamma$ is said to be of {\it diagonal type} if
there exists an orthonormal $\Z$-basis $\{e_1,\dots,e_n\}$ of the lattice $\Lambda$
such that for any element $BL_b\in\Gamma$, $Be_i=\pm e_i$ for $1\le i\le n$ (see \cite{MR3}, Definition 1.3).
Similarly, $M_\G$ is said to be of {\it diagonal type}, if $\G$ is so.
We note that it may be assumed that
the lattice $\Lambda$ of $\Gamma$ is the canonical lattice.

These manifolds have, in particular, holonomy group $F\simeq
\Z_2^r$, for some $1\le r \le n-1$. After conjugation by a
translation we may assume furthermore that  $b\in \frac 12 \Ld$,
for any $BL_b \in \G$ (see \cite{MR3}, Lemma 1.4). In this case,
the terms $e_{\mu,\gamma}$ in the multiplicity formula
(\ref{multip}) become sums of $1$'s and $-1$'s. Moreover, the
traces $\text{tr}_p(B)$ are given by integral values of the so
called Krawtchouk polynomials $K_p^n(x)$ (see \cite{MR2}, Remark
3.6, and also \cite{MR3}; see \cite{KL} for more information on
Krawtchouk polynomials). Namely, we have:
\begin{equation}                                \label{eq.krawtch}
    \text{tr}_p(B)=K_p^n(n-n_B), \quad
    \text{ where } K_p^n(x):= \sum_{t=0}^p (-1)^t \binom xt \binom {n-x}{p-t}.
\end{equation}
For further use, we now give a table with the values of the Krawtchouk polynomials $K_p^n(j)$ for $0\leq p,j \leq 4$.
\begin{equation}\label{krawtvalues}
\begin{tabular}{|c|c|c|c|c|c|}  \hline & & & & &  \\
$p$ & 0 & 1 & 2 & 3 & 4  \\ & & & & & \\  \hline & & & & & \\
$K_p^4(0)$ & 1 & 4 & 6 & 4 & 1 \\ & & & & & \\ \hline
                    & & & & & \\
$K_p^4(1)$ & 1 & 2 & 0 & -2 & -1 \\ & & & & & \\ \hline
                    & & & & & \\
$K_p^4(2)$ & 1 & 0 & -2 & 0 & 1 \\ & & & & & \\ \hline
                    & & & & & \\
$K_p^4(3)$ & 1 & -2 & 0 & 2 & -1 \\ & & & & & \\ \hline
                    & & & & & \\
$K_p^4(4)$ & 1 & -4 & 6 & -4 & 1 \\ & & & & & \\ \hline
\end{tabular}
\end{equation}
\vs

The next lemma gives some
auxiliary facts on the Krawtchouk polynomials.

\begin{lema}[See \cite{KL}, p.\@ 76]\;\;
If $1\le k,j \le n$ we have
\begin{enumerate}
\item[(i)] $ K_k^n(j) = (-1)^j  K_{n-k}^n (j)= (-1)^k  K_k^n (n-j).$
%\item[(ii)] $K_k^n(j) = (-1)^k  K_k^n (n-j)$.
Hence  if $n$ even, then $K_{\frac n2}^n(j)=0$ for $j$ odd and $K_k^n (\frac n2)=0$, for  $k$ odd.
\item[(ii)] $\binom nj K_k^n(j) = \binom nk K_j^n(k)$.
This implies that  $K_k^n(j)=0$ if and only if $K_j^n(k)=0$.
\end{enumerate}
\end{lema}

If $\Gamma$ is of diagonal type and $BL_b \in \G$, we have $n_B=
\text{dim}(\R^n)^B=|\{ 1\le i \le n: Be_i=e_i \}|$. We set
\begin{equation}
 n_B(\tfrac 12) := |\{1 \le i \le n: Be_i=e_i \text { and } b\cdot e_i\equiv \tfrac 12 \mod \Z\}|.
\end{equation}

If $0\le t \le d \le n$, the {\em Sunada numbers} for $\G$  are defined
by
\begin{equation}\label{sunada}
c_{d,t}(\G) :=\big|\big\{BL_b \in \G : n_B=d \text{ and } n_B(\tfrac 12) =t \big \}\big|.
\end{equation}
We note that by Lemma 1.1 in \cite{MR4}, one has that
$n_B(\frac 12)>0$ for any $\g=BL_b \in \G$, $\g \not \in \Ld$.

\begin{rem} In \cite{MR3} it is shown that the equality of the Sunada numbers $c_{d,t}(\G)=c_{d,t}(\G')$  for every $d,t$,  is equivalent to have that $M_\G$ and $M_{\G'}$ verify the conditions in Sunada's theorem. In this case one says that $M_\G$ and $M_{\G'}$
are {\it Sunada isospectral} (see \cite{MR3}, Definition 3.2, Theorem 3.3, and the discussion following it).
In particular this implies that $M_\G$ and $M_{\G'}$ are $p$-isospectral for all~$0\le p\le n$.
\end{rem}

For $\G$ a Bieberbach group, $\tau$ a finite dimensional representation of $\text{O}(n)$ and $Re(s) > 0$, we consider the heat trace zeta function
\begin{equation}                                \label{eq.zetafunc}
    Z^{\Gamma}_\tau(s):=\sum_{\lambda \in \text{$Spec_\tau(M)$}} e^{-\lambda s}
 =\sum_{\mu \ge 0}\, d_{\tau,\mu}(\Gamma)\,
    e^{-4\pi^2 \mu s}                     .
\end{equation}

This series is uniformly convergent for $Re (s) > \varepsilon$, for any $\varepsilon> 0$. In the case when $\tau = \tau_p$, for some $0\le p \le n$,  we write $Z^{\Gamma}_p(s)= Z^{\Gamma}_{\tau_p}(s)$.

%The torsion-free condition on $\G$ implies that for any $BL_b \in \G$ we have $n_B>0$
%= \text{dim}\,\text{ker}(B-\I)>0$
 The lengths of closed geodesics in $M_\G$ are the numbers $\| b_+ + \ld_+ \|$ with $BL_b    \in \G$ running through a full set of representatives of $F= \Lambda\backslash\Gamma$ and $\ld \in \Ld$. Recall that $b_+=p_B(b)$, where $p_B$ denotes the orthogonal projection onto $\text{ker}(B-\I)$.

We will make use of the next result from \cite{MR4}.
\begin{teo}                                 \label{t.zetafunc}
 (i) In the notation above, we have
    \begin{equation}\nonumber                           %\label{eq.longzetaf}
         Z^{\,\Gamma}_\tau(s) = \tfrac 1{|F|} \sum_{BL_b \in F}         \frac{\text{tr}\;\tau(B)}{\hbox{vol}({\Lambda ^*}^{\text{B}})}\,
        (4 \pi s)^{-\frac{n_B}2}
        \sum _{\lambda_+ \in p_B(\Ld)}
        e^{- \frac{ {\| \lambda_+ + b_+\|}^2}{4s}}.
    \end{equation}

$\text{Spec}\,(M_\Gamma)$ determines the lengths of closed geodesics of $M_\G$
and the numbers $n_B$; $\text{Spec}_\tau(M_\Gamma)$ (in particular  $\text{Spec}_p(M_\Gamma)$, for any fixed $p\ge 0$)  determines
the spectrum of the torus $T_\Ld=\Ld\backslash\R^n$ and the cardinality of $F$.

\noi (ii) If $\Gamma$ is a Bieberbach group of diagonal type with
$F\simeq \Z_2^r$ then we have:
%\medskip
    \begin{equation}\nonumber                       %   \label{e.diagzetaf}
         Z^{\,\Gamma}_p(s)=
        \f{1}{2^r} \sum _{d=1}^n K_p^n(n-d)\,
        (4 \pi s)^{-\frac{d}2} \sum _{t=0}^d c_{d,t}(\G)\theta_{d,t}
        (\tfrac 1{4s})
    \end{equation}
%\medskip
where, for $Re(s)>0$,
    \begin{equation}    \nonumber
            \theta_{d,t}(s) = \sum
        _{(m_1,\dots,m_d)\in\Z^d}%{\Sb m_j \in \Z\\ 1\le j \le d \endSb}
        e^{-s \left(\sum_{j=1}^t (\frac 12 + m_j)^2 + \sum_{j=t+1}^d        {m_j}^2\right)}.
    \end{equation}
Furthermore, $M_{\Gamma}$ and $M_{\Gamma'}$ are $p$-isospectral if and only if
$$ K_p^n(n-d)\, c_{d,t}(\G)=K_p^n(n-d)\, c_{d,t}(\G')$$
for each $1\le t\le d \le n$.
In particular, if $c_{d,t}(\G)=c_{d,t}(\G')$ for every $d,t$, then
$M_\G$ and $M_{\G'}$ are $p$-isospectral for all $p$.
If, for $p$ fixed, $K_p^n(x)$ has no integral roots and $M_{\Gamma}$ and $M_{\Gamma'}$
are $p$-isospectral then they are Sunada isospectral, hence $q$-isospectral
for all $0\le q \le n$. In particular, for groups of diagonal type, isospectral implies
 Sunada isospectral.

\end{teo}
\begin{rem}                                 \label{r.rem2}
In \cite{MR4} it is shown by asymptotic methods that the formula
in (i) of the theorem implies that two isospectral flat manifolds
must be both orientable or both nonorientable.
 This is not true for $p$-isospectral pairs (see \cite{MR2}).
\end{rem}

%                           SECTION 2
\section{Bieberbach Groups with canonical translation lattice}
    \label{s.sec2}
    In \cite{BBNWZ} a complete list of isomorphism classes of the crystallographic groups of dimension~4 is given. This list contains 4783 groups out of which    74 are Bieberbach groups.
Among these Bieberbach groups there are 54 which allow the canonical lattice as lattice of translations and 20 that do not. %We will give a full set of representatives for the isometry classes of  manifolds associated to the former groups, yielding 79 non isometric compact flat manifolds. %The fundamental groups

In this section we shall briefly explain how to obtain the Bieberbach
groups from the data given in Tables 1C and 2C in \cite{BBNWZ} and how to
present them in the notation given in Section 1.
%We shall first list 54 Bieberbach groups with canonical lattice of translations and add 25 new ones to account for the different isometry classes. In this way we will obtain $79$  that will give the 79 non-isometric compact flat manifolds.
These will be given in tables containing the matrices $B_i$ and
the corresponding translation vectors $b_i$ and some additional
information. For convenience, in the last table we give the Sunada
numbers for those groups $\G$ having diagonal holonomy
representation.

For simplicity, we shall denote by ${\bf 1},\dots,{\bf 74}$, the 74
Bieberbach groups in Table 1C, p.\@ 408 of \cite{BBNWZ}, giving a full set
of representatives for the isomorphism classes. Groups having the same
holonomy representation (i.e.\@ the same $\Z$-class) are put together into
families ${\mathcal F}_1, \dots, {\mathcal F}_8$. In the notation in
\cite{BBNWZ}, the $\Z$-class is indicated by the first 3 numbers, for
instance, in the group ${\bf 2}$=$2/1/1/2$ below, by $2/1/1$. The full
list is the following: ${\bf 1}$=$1/1/1/1$; ${\bf 2}$=$2/1/1/2$; ${\bf
3}$=$2/1/2/2$; ${\bf 4}$=$2/2/1/2$; ${\bf 5}$=$3/1/1/2$; ${\bf
6}$=$3/1/2/2$; $\!{\mathcal F}_1$=$\{{\bf 7}, {\bf 8},{\bf 9},{\bf
10},{\bf 11}\}$=$\{4/1/1/6, 7, 10, 11, 13\}$; ${\bf 12}$ = $4/1/2/4$;
${\mathcal F}_2$ = $\{{\bf 13},{\bf 14},{\bf 15}\}$ = $\{4/1/3/4, 11,
12\}$; ${\bf 16}$ = $4/1/4/5$; ${\bf 17}$ = $4/1/6/4$; ${\mathcal
F}_3=\{{\bf 18}, {\bf 19}, {\bf 20}, {\bf 21}\}=\{4/2/1/8, 11, 12, 16\}$;
${\bf 22}$=$4/2/3/4$; ${\bf 23}$ = $4/3/1/6$; ${\mathcal F}_4$ = $\{{\bf
24}, {\bf 25},{\bf 26},{\bf 27}\}=\{5/1/2/7, 8, 9, 10\}$; ${\bf 28}$ =
$5/1/3/6$; ${\bf 29}$ = $5/1/4/6$; ${\bf 30}$ = $5/1/6/6$; ${\bf 31}$ =
$5/1/7/4$; ${\bf 32}$ = $5/1/10/4$; ${\mathcal F}_5$ = $\{{\bf 33}, {\bf
34},{\bf 35},{\bf 36},{\bf 37},{\bf 38},{\bf 39},{\bf 40},{\bf 41},{\bf
42}\}$=$\{6/1/1/41, 45, 49, 63,$
 $64, 66, 81, 82, 83, 92\}$;
${\mathcal F}_6$ = $\{{\bf 43},{\bf 44}\}$=$\{6/2/1/27, 50\}$;
${\bf 45}$=$7/2/1/2$; ${\bf 46}$=$7/2/2/2$; ${\bf 47}$=$8/1/1/2$;
${\bf 48}$=$8/1/2/2$; ${\bf 49}$=$9/1/1/2$; ${\bf 50}$=$12/1/2/2$;
${\bf 51}$=$12/1/3/2$; ${\bf 52}$=$12/1/4/2$; ${\bf
53}$=$12/1/6/2$; ${\bf 54}=12/3/4/6$; ${\bf 55}$=$12/3/10/5$;
${\bf 56}$=$12/4/3/11$; ${\mathcal F}_7$=$\{{\bf 57}, {\bf
58}\}$=$\{13/1/1/8, 11\}$; ${\bf 59}$=$13/1/3/8$; ${\mathcal
F}_8$=$\{{\bf 60},{\bf 61},{\bf 62}\}$=$\{13/4/1/14, 20,$ $ 23\}$;
${\bf 63}$=$13/4/4/11$; ${\bf 64}$=$14/1/1/2$; ${\bf
65}$=$14/1/3/2$; ${\bf 66}$=$14/2/3/2$; ${\bf 67}$=$14/3/1/4$;
${\bf 68}$=$14/3/5/4$; ${\bf 69}$=$14/3/6/4$; ${\bf
70}$=$15/1/1/10$; ${\bf 71}$=$15/4/1/10$; ${\bf 72}$=$24/1/2/4$;
 ${\bf 73}$=$24/1/4/4$; ${\bf 74}$=$25/1/1/10$.

\smallskip

In Table 1C in \cite{BBNWZ} we find, for each group $\G$, matrices
$A_1, \dots, A_k$ and translation vectors $a_1, \dots, a_k$
together with a (Bravais) type of lattice. The matrices are not
orthogonal in general, however they are orthogonal with respect to
some suitable inner product. In Table 2C in \cite{BBNWZ}, there is
a list of real symmetric matrices $G= (g_{ij})$ ({\em
Gram~matrices}) corresponding to the types of lattices given in
Table 1C. The entries of these matrices $G$ are parameters
representing the scalar products of the $\Z$-basis vectors of a
lattice of translations $\Ld$ of any group having this lattice
type.

More precisely, Table 1C gives a crystallographic group $\G_1=
\langle A_1L_{a_1},\dots,A_kL_{a_k}, L_{\Ld_1}
 \rangle$ where $\Ld_1=\Z v_1 \oplus \cdots \oplus \Z v_4$
 is a lattice in $\R^4$ %with basis $\{v_1,\dots,v_4\}$.
with $A_i \Ld_1=\Ld_1$ for $1\le i\le 4$. As mentioned before,
Table 2C gives a symmetric matrix $G=(g_{ij})$, where
$g_{ij}=w_i\cdot w_j$ and $\Ld_2=\Z w_1 \oplus \cdots \oplus \Z
w_4$, that gives rise to a group $\G_2= X\G_1 X^{-1}$ for some
unimodular matrix $X$. Such $X\in SL_4(\Z)$ and its inverse
$X^{-1}$ are given in Table 2C, together with the matrix $G$. Now,
the matrices $B_1:=XA_1X^{-1},\dots ,B_k:=XA_kX^{-1}$ obtained in
this way become  orthogonal with respect to the basis $\{w_1,\dots
,w_4\}$ of $\Ld_2$. We will take $\Lambda_2 = \Z^4$. We list the
matrices $X\ne \I$ at the end of the section.

We remark that there is often a difference between the vectors
$a_j$ in Table 1C and those used in this paper. This is so
because, as explained in Section 1, any element $\g \in \G$
decomposes uniquely as $\g=BL_b$, with $B\in \text{O}(n)$ and
$b\in \R^n$, and so $\g$ acts on $x \in \R^n$ by $\g(x)=Bx+Bb$.
However, in \cite{BBNWZ} the elements of $\G$ are written in the
form $\g=(A,a)$ with action $\g(x)=Ax+a$, so $(A,a)$ corresponds
to $AL_{A^{-1}a}$ in our notation. We still have to conjugate
$(A,a)$ by the corresponding matrix $X$. In this way, from the
original data $(A,a)$ we obtain $BL_b = XAX^{-1}L_{XA^{-1}a}$.

\sk In this paper we are interested in groups having the canonical
lattice as lattice of translations. Therefore, we need to look for
the types of lattices having Gram matrices $G$ whose parameters
admit the identity matrix $G=\I$. For example, for the lattices of
type  X/IV (see \cite{BBNWZ}, p.\@ 275) we have
$$G=   \begin{pmatrix}b & 0 & c & c \\
                 & a & 0  & 0 \\
                 &   &  b & -b-2c \\
                 &   &    & b    \end{pmatrix}$$
and we see that $G\ne \I$ for any possible choice of $a,b,c$.
Thus, a  Bieberbach group having type X/IV does not admit the
canonical lattice as translation lattice, and  hence the groups
{\bf 52, 59, 63} are not in the class considered.

Taking into account the full list of Gram matrices in Table 2C in
\cite{BBNWZ} one can check that the 4-dimensional Bieberbach
groups admitting the canonical lattice as lattice of translations
are: {\bf 1}, {\bf 2}, {\bf 3}, {\bf 4}, {\bf 5}, {\bf 6},
${\mathcal F}_1 = \{{\bf 7},{\bf 8},{\bf 9},{\bf 10},{\bf 11}\}$,
{\bf 12}, ${\mathcal F}_2 = \{{\bf 13},{\bf 14},{\bf 15}\}$,
${\mathcal F}_3$ = $\{{\bf 18},{\bf 19},{\bf 20},{\bf 21}\}$, {\bf
22}, {\bf 23}, ${\mathcal F}_4 = \{{\bf 24},{\bf 25},{\bf 26},{\bf
27}\}$, {\bf 28}, {\bf 29},   ${\mathcal F}_5 = \{{\bf 33},{\bf
34},{\bf 35},{\bf 36},{\bf 37}$,${\bf 38},{\bf 39}, {\bf 40},{\bf
41},{\bf 42}\}$, ${\mathcal F}_6 = \{{\bf 43},{\bf 44}\}$, {\bf
45}, {\bf 47}, {\bf 50}, {\bf 51}, {\bf 54}, {\bf 56}, ${\mathcal
F}_7 = \{{\bf 57},{\bf 58}\}$, \nopagebreak[4] ${\mathcal F}_8 =
\{{\bf 60},{\bf 61},{\bf 62}\}$, {\bf 64}, {\bf 67}, {\bf 72} and
{\bf 74}.

We shall present the above 54 groups in tables, together with 25
additional groups to give a complete set of representatives of the
isometry classes of groups in our class. Each one of the  groups
in question lies in one table and groups having the same holonomy
representation are put together in the same table.  Also, for
convenience, all manifolds having holonomy group isomorphic to
$\Z_2$ are put together in the same table (Table 1). Finally, the
groups numbered ${\bf 72}$ and ${\bf 74}$ are not included in the
tables, for simplicity. They can not come up in any
$p$-isospectral pair because they are the only groups with
holonomy groups of order 12 and 24 respectively. The tables show
the non trivial elements in the holonomy group, together with the
corresponding translation vectors and some additional information
that will allow to obtain explicit expressions for the zeta
functions, and to determine the $L$-spectra of the associated
manifolds.

The tables are organized in the following manner. Suppose that $\G
= \langle B_1L_{b_1},\dots,B_kL_{b_k}, \Ld \rangle$. In the first
row of the table, labeled $B$, we give the non trivial matrix
elements of the holonomy group $F$. These matrices $B\in
\text{O}(4,\R)\cap \text{GL}_4(\Z)$ are written in the form
$d(C_1,C_2)$, $d(C_1,C_2,C_3)$ or  $d(C_1,C_2,C_3,C_4)$ where, for
instance, $d(C_1,C_2)$ indicates a  matrix with $C_1, C_2$
 lying diagonally in $B$ and similarly in the other cases.
 Here  $C_i$, $1\le i \le 4$, will be one of
the  following matrices: $\pm 1$, $\pm I$, $\pm \tilde I$, $\pm
J$, $\pm \tilde{J}$, $\pm T$, $\pm T^t$ or $K$ with
$$ I=\begin{pmatrix} 1 & 0 \\ 0 & 1 \end{pmatrix}, \qquad
  \tilde I=\begin{pmatrix} -1 & 0 \\ 0 & 1 \end{pmatrix} \qquad
J=\begin{pmatrix} 0 & 1 \\ 1 & 0 \end{pmatrix}, \qquad
 \tilde J=\begin{pmatrix} 0 & 1 \\ -1 & 0 \end{pmatrix}, $$
$$    T=\begin{pmatrix} 0 & 1 & 0 \\
                 0 & 0 & 1 \\

                 1 & 0 & 0 \end{pmatrix}, \qquad
    K=\begin{pmatrix} 0 & 0 & 1
             \\ 0 & 1 & 0 \\
                1 & 0 & 0 \end{pmatrix}.$$

The second group of entries, corresponds to the translation
vectors $b$. Here we include one row for each group in the family,
denoting with primes the different isometry classes in the same
diffeomorphism class.

In the table we also indicate a set of  generators over $\Z$ of
the space fixed by $B$, i.e.\@ $\Ld^B=\{ v\in \Ld :  Bv=v\}$, the
numbers $n_B$, $vol(\Ld^B)$ and the $p$-traces $tr_p(B)$.  In
general, these traces are given by integer values of Krawtchouk
polynomials, $K_p^n(j)$, but there are some exceptions. In these
cases all the values of the $p$-traces,  for $0\leq p\leq 4$, are
listed. The components of the vectors $b$ that are not in $\Ld^B$
are put in between parentheses in the tables. In the next section
we will use these tables to compute, for each group $\G$, the
corresponding $p$-heat trace zeta function from which one can read
off spectral information on $M_\G$.

Given a group ${\bf k}$, there may be several different isometry
classes of Bieberbach groups in the isomorphism  class of ${\bf
k}$, to be denoted by ${\bf k'}$, ${\bf k''}$ and so on (see for
instance, Table 2.1). One shows that all such classes can be
obtained by conjugation of the matrices in the point group $F$ of
$k$, by some $C\in \text{GL}_4(\Z)\cap Z(F) \smallsetminus
\text{O}(4)$, where $Z(F)$ denotes the centralizer of $F$.

The following is a list of representatives for the isometry
classes of Bieberbach groups having canonical lattice of
translations. The corresponding classes in the notation in
\cite{BBNWZ} are given in between parentheses.

\medskip

 \noindent
{\sl Table 2.1}: ${\bf 2}\, (2/1/1), {\bf 3}\,(2/1/2), {\bf 4}\,(2/2/1), {\bf 5}\,(3/1/1), {\bf 6}\, (3/1/2)$  $(F\simeq \Z_2)$.
\

\begin{center}
 \begin{tabular}{|c||*{5}{c|}}
\hline & & & & & \\
  \#    & \bf2 & \bf3 & \bf4    & \bf5 & \bf6                  \\  & & & & & \\             \hline & & & & & \\
  $B$ & $d(I,1,-1)$ & $d(I,J)$ & $d(-I,-1,1)$ & $d(-I,I)$ & $d(1,J,-1)$ \\                                     & & & & & \\
            \hline & & & & & \\
  $b$ & $\frac{e_3}{2}$ & $\f{e_1}{2}$ & $\f{e_4}{2}$ & $\f{e_4}{2}$ & $\f{e_1}{2}$                         \\ & & & & & \\
  $b'$ & $\frac{e_2+e_3}{2}$ & $\f{e_1+e_2}{2}$ &  & $\f{e_3+e_4}{2}$ & $\f{e_1+e_2+e_3}{2}$        \\ & & & & & \\
  $b''$ & $\frac{e_1+e_2+e_3}{2}$ & $\f{e_1+e_3+e_4}{2}$ &  &  &        \\ & & & & & \\
  $b'''$ &  & $\f{e_1+e_2+e_3+e_4}{2}$ &  &  &      \\ & & & & & \\
                          \hline & & & & & \\
  $\Lambda^B$ & $e_1, e_2, e_3$ & $e_1, e_2, e_3 + e_4$ & $e_4$ & $e_3, e_4$ & $e_1, e_2 + e_3$                                         \\ & & &  & & \\
  $n_B$ & 3 & 3 & 1 & 2 & 2                             \\ & & &  & & \\
  $vol(\Lambda^B)$ & 1 & $\sqrt{2}$ & 1 & 1 & $\sqrt{2}$        \\ & & &  & & \\
                                             \hline  & & &  & & \\
$tr_p(B)$ &   $K_p^4(1)$  & $K_p^4(1)$  & $K_p^4(3)$  & $K_p^4(2)$  & $K_p^4(2)$                                                                                                    \\ & & &  & & \\
  \hline    & & & & &       \\
$\beta_1, \beta_2$ & $3,3$ & $3,3$ &  $1,3$ &  $2,2$ & $2,2$
 \\ & & & & & \\
{\em Orientable}    & {\em no} & {\em no} & {\em no} & {\em yes} & {\em yes} \\ & & & & & \\
{\em Diagonal-type} & {\em yes} & {\em no} & {\em yes} & {\em yes} & {\em no} \\ & & & & & \\ \hline
 \end{tabular}
\end{center}
\medskip

\medskip \noindent
{\sl Table 2.2}: {\bf 47}\,(8/1/1) ($F\simeq \Z_3$, {\em non-diagonal, orientable}, $\beta_1 =2, \beta_2=1$)
\

\begin{center}
 \begin{tabular}{|c||c|c|}
\hline & & \\
& $\g_1$ & $\g_1^2$\\ && \\\hline & & \\
  $B$ & $d(1,T)$ & $d(1,T^t)$                           \\ & &  \\ \hline & &  \\
  $b_{\bf 47}$ & $\frac{e_1}{3}$ & $2\cdot \f{e_1}{3}$      \\ & &  \\  $b_{\bf 47'}$ & $\frac{e_1+e_2+e_3+e_4}{3}$ & $2\cdot\f{e_1+e_2+e_3+e_4}{3}$        \\ & &  \\  \hline & &  \\
  $\Lambda^B$ & $e_1, e_2+ e_3+e_4$ & $e_1, e_2+ e_3 + e_4$  \\  & & \\
  $n_B$ & 2 & 2                                     \\  & & \\
  $vol(\Lambda^B)$ & $\sqrt{3}$ & $\sqrt{3}$        \\  & & \\  \hline  & & \\
  $tr_p(B)$ & 1 1 0 1 1 & 1 1 0 1 1                 \\ & &  \\
  \hline
 \end{tabular}

\end{center}
\pagebreak

\bigskip

\noindent
{\sl Table 2.3}: ${\mathcal F}_1=\{ {\bf 7}, {\bf 8}, {\bf 9}, {\bf 10}, {\bf 11} \}\, (4/1/1)$    ($F\simeq \Z_2^2$, {\em diagonal, non-orientable}, $\beta_1=2, \beta_2=1$)
\begin{center}
 \begin{tabular}{|c||c|c|c|}
    \hline & & & \\ & $\g_1$ & $\g_2$ & $\g_1\g_2$ \\ & & & \\  \hline & & & \\
  $B$ & $d(I,1,-1)$ & $d(I,-1,1)$ & $d(I,-I)$      \\ & & & \\ \hline  & & & \\
    $b_{\bf 7}$ & $\frac{e_3}{2}$ & $\f{e_2}{2}$ & $\f{e_2 + (e_3)}{2}$ \\ & & &  \\
    $b_{\bf 7'}$ & $\frac{e_3}{2}$ & $\f{e_1+e_2}{2}$ & $\f{e_1+e_2 + (e_3)}{2}$ \\ & & &  \\
    $b_{\bf 8}$ & $\frac{e_3}{2}$ & $\f{e_2+e_4}{2}$ & $\f{e_2 + (e_3+e_4)}{2}$ \\ & & &  \\
    $b_{\bf 8'}$ & $\frac{e_3}{2}$ & $\f{e_1+e_2+e_4}{2}$ & $\f{e_1+e_2 + (e_3+e_4)}{2}$ \\ & & &  \\
    $b_{\bf 9}$ & $\frac{e_2}{2}$ & $\f{e_1}{2}$ & $\f{e_1 +e_2}{2}$ \\ & & &  \\
      $b_{\bf 9'}$ & $\frac{e_2}{2}$ & $\f{e_1+e_2}{2}$ & $\f{e_1}{2}$ \\ & & &  \\
    $b_{\bf 10}$ & $\frac{e_2}{2}$ & $\f{e_1+e_4}{2}$ & $\f{e_1+e_2+(e_4)}{2}$ \\ & & &  \\
    $b_{\bf 10'}$ & $\frac{e_2}{2}$ & $\f{e_1+e_2+e_4}{2}$ & $\f{e_1+(e_4)}{2}$ \\ & & &  \\
    $b_{\bf 10''}$ & $\frac{e_1+e_2}{2}$ & $\f{e_1+e_4}{2}$ & $\f{e_2+(e_4)}{2}$ \\ & & &  \\
    $b_{\bf 11}$ & $\frac{e_2+e_3}{2}$ & $\f{e_1+e_4}{2}$ & $\f{e_1+e_2 + (e_3+e_4)}{2}$ \\ & & &  \\
    $b_{\bf 11'}$ & $\frac{e_2+e_3}{2}$ & $\f{e_1+e_2+e_4}{2}$ & $\f{e_1+(e_3+e_4)}{2}$ \\ & & &  \\                              \hline & & & \\
  $\Lambda^B$ & $e_1, e_2, e_3$ & $e_1, e_2, e_4$ & $e_1, e_2$    \\ & & & \\
  $n_B$ & 3 & 3 & 2                                          \\ & & & \\
  $vol(\Lambda^B)$ & 1 & 1 & 1                      \\ & & & \\ \hline  & & & \\
  $tr_p(B)$ & $K_p^4(1)$  & $K_p^4(1)$  & $K_p^4(2)$             \\ & & &  \\
  \hline
 \end{tabular}

\end{center}

\medskip \noindent
{\sl Table 2.4}: {\bf 12}\,(4/1/2) ($F\simeq \Z_2^2$, {\em non-diagonal, non-orientable}, $\beta_1=2, \beta_2=1$)

\begin{center}
 \begin{tabular}{|c||c|c|c|}
\hline & & & \\& $\g_1$ & $\g_2$ & $\g_1\g_2$ \\ & & & \\  \hline & & & \\
  $B$ & $d(I,J)$ & $d(I,-J)$ & $d(I,-I)$                    \\ & & &  \\ \hline  & & &  \\
  $b_{\bf 12}$ & $\frac{e_2}{2}$ & $\f{e_1}{2}$ & $\f{e_1 + e_2}{2}$ \\ & & &  \\       $b_{\bf 12'}$ & $\frac{e_1+e_2}{2}$ & $\f{e_1}{2}$ & $\f{e_2}{2}$ \\ & & &  \\                                                  \hline & & & \\
  $\Lambda^B$ & $e_1, e_2, e_3+e_4$ & $e_1, e_2, e_3-e_4$ & $e_1, e_2$    \\ & & & \\
  $n_B$ & 3 & 3 & 2                                          \\ & & & \\
  $vol(\Lambda^B)$ & $\sqrt{2}$ & $\sqrt{2}$ & 1                    \\ & & & \\ \hline  & & & \\
  $tr_p(B)$ & $K_p^4(1)$  & $K_p^4(1)$  & $K_p^4(2)$             \\ & & &  \\
  \hline
 \end{tabular}
\end{center}
\pagebreak

\noindent {\sl Table 2.5}: ${\mathcal F}_2  =\{ {\bf 13}, {\bf
14}, {\bf 15} \}\, (4/1/3)$  ($F\simeq  \Z_2^2$, {\em
non-diagonal, non-orientable}, $\beta_1=1, \beta_2=0$)

\begin{center}
 \begin{tabular}{|c||c|c|c|}
  \hline  & & & \\& $\g_1$ & $\g_2$ & $\g_1\g_2$ \\ & & & \\  \hline & & & \\
$B$ & $d(1,I,-1)$ & $d(1,J,1)$ & $d(1,J,-1)$ \\ & & &  \\  \hline & & & \\
$b_{\bf13}$ & $\frac{e_1}2$ & $\f{e_4}{2}$ & $\f{e_1+(e_4)}{2}$
\\ & & & \\
$b_{\bf13'}$ & $\frac{e_1+e_2+e_3}2$ & $\f{e_4}{2}$ & $\f{e_1+e_2+e_3+(e_4)}{2}$
\\ & & & \\
$b_{\bf14}$ & $\frac{e_2+e_3}{2}$ & $\f{e_1}{2}$ & $\f{e_1+e_2+e_3}{2}$
\\ & & & \\
$b_{\bf14'}$ & $\frac{e_2+e_3}{2}$ & $\f{e_1+e_2+e_3}{2}$ & $\f{e_1}{2}$
\\ & & & \\
$b_{\bf15}$ & $\frac{e_2+e_3}{2}$ & $\f{e_1+e_4}{2}$ & $\f{e_1+e_2+e_3+(e_4)}{2}$
\\ & & & \\
$b_{\bf15'}$ & $\frac{e_2+e_3}{2}$ & $\f{e_1+e_2+e_3+e_4}{2}$ & $\f{e_1+(e_4)}{2}$
\\ & & & \\ \hline & & & \\
  $\Lambda^B$ & $e_1, e_2, e_3$ & $e_1, e_2+e_3, e_4$ & $e_1,e_2+e_3$
\\ & & &  \\
  $n_B$ & 3 & 3 & 2 \\ & & & \\
  $vol(\Lambda^B)$ & 1 & $\sqrt{2}$ & $\sqrt{2}$ \\ & & & \\ \hline & & &  \\
  $tr_p(B)$ & $K_p^4(1)$  & $K_p^4(1)$  & $K_p^4(2)$ \\ & & & \\
  \hline
 \end{tabular}
\end{center}

\medskip
\noindent
{\sl Table 2.6}: ${\mathcal F}_3  =\{{\bf 18}, {\bf 19}, {\bf 20}, {\bf 21}\} \,\, (4/2/1)$  ($F\simeq  \Z_2^2$, {\em diagonal, non-orientable}, $\beta_1=1, \beta_2=1$)

\begin{center}
 \begin{tabular}{|c||c|c|c|}  \hline  & & & \\& $\g_1$ & $\g_2$ & $\g_1\g_2$ \\ & & & \\  \hline & & & \\
  $B$ & $d(I,-1,1)$ & $d(-I,-1,1)$ & $d(-I,I)$ \\ & & &  \\ \hline & & & \\
  $b_{\bf18}$ & $\frac{e_4}2$ & $\f{(e_3)+e_4}{2}$ &
 $\f{e_3}{2}$ \\ & & & \\
  $b_{\bf19}$ & $\frac{e_2}{2}$ & $\f{e_4}{2}$ & $\f{(e_2)+e_4}{2}$ \\ & & & \\
  $b_{\bf19'}$ & $\frac{e_1+e_2}{2}$ & $\f{e_4}{2}$ & $\f{(e_1+e_2)+e_4}{2}$ \\ & & & \\
  $b_{\bf20}$ & $\frac{e_2}{2}$ & $\f{(e_3)+e_4}{2}$ & $\f{(e_2)+e_3+e_4}{2}$
\\ & & & \\
$b_{\bf20'}$ & $\frac{e_1+e_2}{2}$ & $\f{(e_3)+e_4}{2}$ & $\f{(e_1+e_2)+e_3+e_4}{2}$
\\ & & & \\
  $b_{\bf21}$ & $\frac{e_2+e_4}{2}$ & $\f{(e_3)+e_4}{2}$ & $\f{(e_2)+e_3}{2}$  \\ & & & \\
 $b_{\bf21'}$ & $\frac{e_1+e_2+e_4}{2}$ & $\f{(e_3)+e_4}{2}$ & $\f{(e_1+e_2)+e_3}{2}$  \\ & & & \\ \hline & & & \\
  $\Lambda^B$ & $e_1, e_2, e_4$ & $e_4$ & $e_3,e_4$ \\ & & &  \\
  $n_B$ & 3 & 1 & 2 \\ & & &  \\
  $vol(\Lambda^B)$ & 1 & 1 & 1 \\ & & & \\ \hline
& & &  \\
  $tr_p(B)$ & $K_p^4(1)$  & $K_p^4(3)$  & $K_p^4(2)$ \\
   & & & \\
  \hline
 \end{tabular}
\end{center}

%\pagebreak
%\medskip

\pagebreak

 \noindent{\sl Table 2.7}: {\bf 22}\, (4/2/3) ($F\simeq \Z_2^2$, {\em non-diagonal, non-orientable}, $\beta_1=1, \beta_2=0$)

\begin{center}
 \begin{tabular}{|c||c|c|c|}
    \hline & & & \\& $\g_1$ & $\g_2$ & $\g_1\g_2$ \\ & & & \\  \hline & & & \\
  $B$ & $d(1,J,1)$ & $d(-1,-I,,1)$ & $d(-1,-J,1)$                   \\ & & & \\ \hline & & & \\
  $b_{\bf 22}$ & $\frac{e_1}{2}$ & $\f{e_4}{2}$ & $\f{(e_1) + e_4}{2}$ \\ & & &  \\      $b_{\bf 22'}$ & $\frac{e_1+e_2+e_3}{2}$ & $\f{e_4}{2}$ & $\f{(e_1+e_2+e_3)+e_4}{2}$ \\ & & &  \\                               \hline & & & \\
  $\Lambda^B$ & $e_1, e_2+e_3, e_4$ & $e_4$ & $e_2-e_3, e_4$      \\ & & & \\
  $n_B$ & 3 & 1 & 2                                          \\ & & & \\
  $vol(\Lambda^B)$ & $\sqrt{2}$ & 1 & $\sqrt{2}$                        \\ & & & \\ \hline  & & & \\
  $tr_p(B)$ & $K_p^4(1)$  & $K_p^4(3)$  & $K_p^4(2)$             \\ & & &  \\
  \hline
 \end{tabular}
\end{center}

\medskip
\noindent
 {\sl Table 2.8}: {\bf 23}\, (4/3/1) ($F\simeq \Z_2^2$, {\em diagonal,
non-orientable}, $\beta_1=0, \beta_2=1$)

\begin{center}
 \begin{tabular}{|c||c|c|c|}
\hline & & & \\& $\g_1$ & $\g_2$ & $\g_1\g_2$ \\ & & & \\  \hline & & & \\
  $B$ & $d(-I,1,-1)$ & $d(-I,-1,1)$ & $d(I,-I)$                     \\ & & & \\ \hline & & & \\  $b_{\bf 23}$ & $\frac{e_3}{2}$ & $\f{(e_2)+e_4}{2}$ & $\f{e_2+(e_3+ e_4)}{2}$ \\ & & &  \\   $b_{\bf 23'}$ & $\frac{e_3}{2}$ & $\f{(e_1+e_2)+e_4}{2}$ & $\f{e_1+e_2+(e_3+ e_4)}{2}$ \\ & & &  \\                               \hline & & & \\
  $\Lambda^B$ & $e_3$ & $e_4$ & $e_1, e_2$    \\ & & & \\
  $n_B$ & 1 & 1 & 2                                          \\ & & & \\
  $vol(\Lambda^B)$ & 1 & 1 & 1                      \\ & & & \\ \hline  & & & \\
  $tr_p(B)$ & $K_p^4(3)$  & $K_p^4(3)$  & $K_p^4(2)$             \\ & & &  \\
  \hline
 \end{tabular}
\end{center}

\medskip
\noindent
{\sl Table 2.9}: ${\mathcal F}_4  =\{{\bf 24}, {\bf25}, {\bf 26}, {\bf 27}\}\, (5/1/2)$ ($F\simeq  \Z_2^2$, {\em diagonal, orientable}, $\beta_1=1, \beta_2=0$)

\begin{center}
 \begin{tabular}{|c||c|c|c|}
  \hline  & & & \\& $\g_1$ & $\g_2$ & $\g_1\g_2$ \\ & & & \\  \hline & & & \\
$B$ & $d(-1,-1,1,1)$ & $d(1,-1,-1,1)$ & $d(-1,1,-1,1)$
\\& & &  \\ \hline & & & \\
$b_{\bf24}$ & $\frac{e_4}2$ & $\f{(e_2)+e_4}{2}$ &
 $\f{e_2}{2}$
\\ & & & \\
$b_{\bf25}$ & $\frac{e_4}{2}$ & $\f{e_1+(e_2)}{2}$ & $\f{(e_1)+e_2+e_4}{2}$ \\ & & & \\
$b_{\bf26}$ & $\frac{e_3}{2}$ & $\f{e_1+(e_2)}{2}$ & $\f{(e_1)+e_2+(e_3)}{2}$
\\ & & & \\
$b_{\bf27}$ & $\frac{e_3}{2}$ & $\f{e_1+(e_2)+e_4}{2}$ & $\f{(e_1)+e_2+(e_3)+e_4}{2}$  \\ & & & \\
\hline  & & & \\
$\Lambda^B$ & $e_3, e_4$ & $e_1, e_4$ & $e_2, e_4$
\\ & & &  \\
$n_B$ & 2 & 2 & 2 \\ & & &  \\
$vol(\Lambda^B)$ & 1 & 1 & 1  \\ & & & \\ \hline
& & &  \\
$tr_p(B)$ & $K_p^4(2)$  & $K_p^4(2)$  & $K_p^4(2)$ \\
& & & \\ \hline
\end{tabular}
\end{center}

\pagebreak
 \noindent
 {\sl Table 2.10}: {\bf 28}\, (5/1/3) ($F\simeq \Z_2^2$, {\em non-diagonal, orientable}, $\beta_1=1, \beta_2=0$)
   \begin{center}
 \begin{tabular}{|c||c|c|c|}                    \hline & & & \\& $\g_1$ & $\g_2$ & $\g_1\g_2$ \\ & & & \\  \hline & & & \\
  $B$ & $d(1,-J,-1)$ & $d(1,J,-1)$ & $d(1,-I,1)$                    \\ & & & \\
  $b_{\bf 28}$ & $\frac{e_1}{2}$ & $\f{e_1 +(e_4)}{2}$ & $\f{e_4}{2}$ \\ & & &  \\                  \hline & & & \\
  $\Lambda^B$ & $e_1, e_2-e_3$ & $e_1, e_2+e_3$ & $e_1, e_4$      \\ & & & \\
  $n_B$ & 2 & 2 & 2                                          \\ & & & \\
  $vol(\Lambda^B)$ & $\sqrt{2}$ & $\sqrt{2}$ & 1                        \\ & & & \\ \hline  & & & \\
  $tr_p(B)$ & $K_p^4(2)$  & $K_p^4(2)$  & $K_p^4(2)$             \\ & & &  \\
  \hline
 \end{tabular}
\end{center}

\medskip \noindent
 {\sl Table 2.11}: {\bf 29}\, (5/1/4) ($F\simeq \Z_2^2$, {\em non-diagonal, orientable}, $\beta_1=1, \beta_2=0$)

  \begin{center}
 \begin{tabular}{|c||c|c|c|}                            \hline & & & \\& $\g_1$ & $\g_2$ & $\g_1\g_2$ \\ & & & \\  \hline & & & \\
  $B$ & $d(-1,J,1)$ & $d(1,J,-1)$ & $d(-1,1,1,-1)$                      \\ & & & \\ \hline & & & \\
  $b_{\bf 29}$ & $\frac{(-e_2+e_3)+e_4}{2}$ & $\f{e_1}{2}$ & $\f{(e_1)+e_2 +e_3+(e_4)}{2}$ \\ & & &  \\
 $b_{\bf 29'}$ & $\frac{e_3+e_4}{2}$ & $\f{e_1}{2}$ & $\f{(e_1)+e_3+(e_4)}{2}$ \\ & & &  \\                                 \hline & & & \\
  $\Lambda^B$ & $e_2+e_3, e_4$ & $e_1, e_2+e_3$ & $e_2, e_3$      \\ & & & \\
  $n_B$ & 2 & 2 & 2                                          \\ & & & \\
  $vol(\Lambda^B)$ & $\sqrt{2}$ & $\sqrt{2}$ & 1                        \\ & & & \\ \hline  & & & \\
  $tr_p(B)$ & $K_p^4(2)$  & $K_p^4(2)$  & $K_p^4(2)$             \\ & & &  \\
  \hline
 \end{tabular}
\end{center}

\medskip
\noindent
 {\sl Table 2.12}: {\bf 45}\, (7/2/1) ($F\simeq \Z_4$, {\em non-diagonal, orientable}, $\beta_1=2, \beta_2=1$)

   \begin{center}
 \begin{tabular}{|c||c|c|c|}                            \hline & & & \\& $\g_1$ & $\g_1^2$ & $\g_1^3$ \\ & & & \\  \hline & & & \\
  $B$ & $d(I,-\tilde J)$ & $d(I,-I)$ & $d(I,\tilde J)$                      \\ & & & \\ \hline & & & \\
  $b_{\bf 45}$ & $\frac{e_2}{4}$ & $\f{e_2}{2}$ & $\f{3e_2}{4}$  \\ & & &  \\    $b_{\bf 45'}$ & $\frac{e_1+e_2}{4}$ & $\f{e_1+e_2}{2}$ & $\f {3e_1+3e_2}4$  \\ & & &  \\                                       \hline & & & \\
  $\Lambda^B$ & $e_1, e_2$ & $e_1, e_2$ & $e_1, e_2$       \\ & & & \\
  $n_B$ & 2 & 2 & 2                                          \\ & & & \\
  $vol(\Lambda^B)$ & $1$ & $1$ & 1                        \\ & & & \\ \hline  & & & \\
  $tr_p(B)$ & 1 2 2 2 1  & $K_p^4(2)$  & 1 2 2 2 1             \\ & & &  \\
  \hline
 \end{tabular}
\end{center}

 \noindent

\pagebreak
 {\sl Table 2.13}: {\bf 50} \, (12/1/2)($F\simeq \Z_4$, {\em non-diagonal, non-orientable}, $\beta_1=1, \beta_2=0$)

   \begin{center}
 \begin{tabular}{|c||c|c|c|}            \hline & & & \\& $\g_1$ & $\g_1^2$ & $\g_1^3$ \\ & & & \\  \hline & & & \\
  $B$ & $d(-1,1,-\tilde J)$ & $d(I,-I)$ & $d(-1,1,\tilde J)$                    \\ & & & \\
  $b_{\bf 50}$ & $\frac{e_2}{4}$ & $\f{e_2}{2}$ & $\f{3e_2}{4}$  \\ & & &  \\                                                   \hline & & & \\
  $\Lambda^B$ & $e_2$ & $e_1, e_2$ & $e_2$    \\ & & & \\
  $n_B$ & 1 & 2 & 1                                          \\ & & & \\
  $vol(\Lambda^B)$ & 1 & 1 & 1                      \\ & & & \\ \hline  & & & \\
  $tr_p(B)$ & 1 0 0 0 -1  & $K_p^4(2)$  & 1 0 0 0 -1             \\ & & &  \\
  \hline
 \end{tabular}
\end{center}

\medskip \noindent
 {\sl Table 2.14}: {\bf 51} \, (12/1/3)($F\simeq \Z_4$, {\em non-diagonal, non-orientable}, $\beta_1=0, \beta_2=1$)
   \begin{center}
 \begin{tabular}{|c||c|c|c|}                            \hline & & & \\& $\g_1$ & $\g_1^2$ & $\g_1^3$ \\ & & & \\  \hline & & & \\
  $B$ & $d(J,-\tilde J)$ & $d(I,-I)$ & $d(J,\tilde J)$                      \\ & & & \\
  $b_{\bf 51}$ & $\frac{e_2}{2}$ & $\f{e_1+e_2}{2}$ & $\f{e_1}{2}$  \\ & & &  \\                    \hline & & & \\
  $\Lambda^B$ & $e_1+e_2$ & $e_1, e_2$ & $e_1+e_2$    \\ & & & \\
  $n_B$ & 1 & 2 & 1                                          \\ & & & \\
  $vol(\Lambda^B)$ & $\sqrt{2}$ & 1 & $\sqrt{2}$                        \\ & & & \\ \hline  & & & \\
  $tr_p(B)$ & 1 0 0 0 -1  & $K_p^4(2)$  & 1 0 0 0 -1             \\ & & &  \\
  \hline
 \end{tabular}
\end{center}

\medskip \noindent
 {\sl Table 2.15}: {\bf 64} \, (14/1/1)($F\simeq \Z_6$, {\em non-diagonal, non-orientable}, $\beta_1=1, \beta_2=0$)
   \begin{center}
 \begin{tabular}{|c||c|c|c|c|c|}                            \hline & & & & &\\& $\g_1$ & $\g_1^2$ & $\g_1^3$ & $\g_1^4$ & $\g_1^5$ \\ & & & & & \\  \hline  & & & & & \\
  $B$ & $d(1,-T)$ & $d(1,T^t)$ & $d(1,-1,-1,-1)$ & $d(1,T)$ & $d(1,-T^t)$                    \\ & & & & & \\
  $b_{\bf 64}$ & $\frac{e_1}{6}$ & $\f{e_1}{3}$ & $\f{e_1}{2}$  & $\f{2e_1}{3}$ & $\f{5e_1}{6}$ \\ & & & & & \\                                     \hline & & & & & \\
  $\Lambda^B$ & $e_1$ & $e_1, e_2\!+\!e_3\!+\!e_4$ & $e_1$ & $e_1, e_2\!+\!e_3\!+\!e_4$ &   $e_1$   \\ & & & & &  \\
  $n_B$ & 1 & 2 & 1   & 2 & 1                                        \\ & & &   & & \\
  $vol(\Lambda^B)$ & 1 & $\sqrt{3}$ & 1 & $\sqrt{3}$ & 1                        \\ & & & & & \\ \hline  & & & & & \\
  $tr_p(B)$ & 1 1 0  -1 1  & 1 1 0 1 1 &  $K_p^4(3)$ & 1 1 0  1 1   & 1 1 0  -1 1             \\ & & & & &  \\
  \hline
 \end{tabular}
\end{center}

%\pagebreak
\medskip
\noindent
 {\sl Table 2.16}: {\bf 67} \, (14/3/1)($F\simeq D_3$, {\em non-diagonal, orientable}, $\beta_1=0, \beta_2=0$)
   \begin{center}
 \begin{tabular}{|c||c|c|c|c|c|}                            \hline & & & & &\\& $\g_1$ & $\g_1^2$ & $\g_2$ & $\g_1 \g_2$ & $\g_1^2 \g_2$ \\ & & & & & \\  \hline  & & & & & \\
  $B$ & $d(1,T)$ & $d(1,T^t)$ & $d(-1,J,1)$ & $d(-1,1,J)$ & $d(-1,K)$                   \\ & & & & & \\
  $b_{\bf 67}$ & $\frac{e_1}{3}\!+\!\f{e_3+e_4}{2}$ & $\f{2e_1}{3}\!+\!\f{e_2+e_3}{2}$ & $\f{e_4}{2}$  & $(\!\f{e_1}{3}\!)\!+\!\f{e_2}{2}$ & $(\!\f{2e_1}{3}\!)\!+\!\f{e_2+e_3+e_4}{2}$ \\ & & & & & \\                                                     \hline & & & & & \\
  $\Lambda^B$ & $e_1, e_2\!\!+\!\!e_3\!\!+\!\!e_4$ & $e_1, e_2\!\!+\!\!e_3\!\!+\!\!e_4$ & $e_2\!+\!e_3, e_4$ & $e_2, e_3\!+\!e_4$ &     $e_3, e_2\!+\!e_4$   \\ & & & & &  \\
  $n_B$ & 2 & 2 & 2   & 2 & 2                                        \\ & & &   & & \\
  $vol(\Lambda^B)$ & $\sqrt{3}$  & $\sqrt{3}$ & $\sqrt{2}$  & $\sqrt{2}$ & $\sqrt{2}$                       \\ & & & & & \\ \hline  & & & & & \\
  $tr_p(B)$ & 1 1 0 1 1  & 1 1 0 1 1 &  $K_p^4(2)$ & $K_p^4(2)$    & $K_p^4(2)$              \\ & & & & &  \\
  \hline
 \end{tabular}
\end{center}

\pagebreak
 \noindent
 {\sl Table 2.17}: ${\mathcal F}_5  =\{{\bf
33}, {\bf 34}, {\bf 35}, {\bf 36}, {\bf 37}, {\bf 38}, {\bf 39},
{\bf 40}, {\bf 41}, {\bf 42}\}\, (6/1/1)$  ($F\simeq  \Z_2^3$,
{\em diagonal, non-orient\-able}, $\beta_1=1, \beta_2=0$).

\begin{center}  \scriptsize{
 \begin{tabular}{|c||c|c|c|c|c|c|c|}
        \hline  & & & & & & & \\& $\g_1$ & $\g_2$ & $\g_3$ & $\g_1\g_2$ & $\g_1 \g_3$ & $\g_2 \g_3$ & $\g_1 \g_2 \g_3$ \\& & & & & & & \\  \hline & & & & & & & \\
  $B$  &  $d(\tilde I,I)$  &  $d(I,\tilde I)$  &  $d(-I,I)$ &  $d(\tilde I, \tilde I)$  &  $d(-\tilde I,I)$  &  $d(-I,\tilde I)$ & $d(-\tilde I, \tilde I)$ \\    & & & & & & & \\ \hline & & & & & & & \\
 $b_{\bf33}$     & $\frac{e_4}{2}$ & $\f{e_2}{2}$ &
           $\f{(e_1)+e_4}{2}$     & $\frac{e_2+e_4}{2}$ &
           $\f{e_1}{2}$ & $\f{(e_1+e_2)+e_4}{2}$  &
           $\frac{e_1+(e_2)}{2}$      \\ & & & & & & & \\
 $b_{\bf34}$        & $\frac{e_4}{2}$ & $\f{e_2}{2}$ &
        $\f{(e_1)+e_3+e_4}{2}$     & $\frac{e_2+e_4}{2}$ &
         $\f{e_1+e_3}{2}$ & $\f{(e_1+e_2+e_3)+e_4}{2}$ &
         $\frac{e_1+(e_2+e_3)}{2}$   \\ & & & & & & & \\
$b_{\bf35}$     & $\frac{e_4}{2}$ & $\f{e_2+e_4}{2}$ &
           $\f{(e_1)+e_3}{2}$     & $\frac{e_2}{2}$ &
           $\f{e_1+e_3+e_4}{2}$ &
            $\f{(e_1+e_2+e_3)+e_4}{2}$  &
           $\frac{e_1+(e_2+e_3)}{2}$  \\ & & & & & & & \\
 $b_{\bf36}$        & $\frac{e_3}{2}$ & $\f{e_2}{2}$ &
        $\f{e_4}{2}$     & $\frac{e_2+(e_3)}{2}$ &
         $\f{e_3+e_4}{2}$ & $\f{(e_2)+e_4}{2}$ &
         $\frac{(e_2+e_3)+e_4}{2}$   \\ & & & & & & & \\
$b_{\bf37}$     & $\frac{e_3}{2}$ & $\f{e_2}{2}$ &
           $\f{(e_1)+e_4}{2}$     & $\frac{e_2+(e_3)}{2}$ &
          $\f{e_1+e_3+e_4}{2}$ & $\f{(e_1+e_2)+e_4}{2}$  &
        $\frac{e_1+(e_2+e_3)+e_4}{2}$  \\& & & &  & & &\\
 $b_{\bf38}$        & $\frac{e_3}{2}$ & $\f{e_2}{2}$ &
        $\f{(e_1)+e_3+e_4}{2}$   & $\frac{e_2+(e_3)}{2}$ &
         $\f{e_1+e_4}{2}$ & $\f{(e_1+e_2+e_3)+e_4}{2}$ &
       $\frac{e_1+(e_2)+e_4}{2}$    \\& & & & &  & &\\
$b_{\bf39}$     & $\frac{e_3}{2}$ & $\f{e_1+e_2}{2}$ &
           $\f{e_4}{2}$     & $\frac{(e_1)+e_2+(e_3)}{2}$ &
           $\f{e_3+e_4}{2}$ & $\f{(e_1+e_2)+e_4}{2}$  &
        $\frac{e_1+(e_2+e_3)+e_4}{2}$ \\& & & & & & &\\
 $b_{\bf40}$       & $\frac{e_3}{2}$ & $\f{e_1+e_2}{2}$ &
        $\f{(e_1)+e_4}{2}$  & $\frac{(e_1)+e_2+(e_3)}{2}$ &
         $\f{e_1+e_3+e_4}{2}$ & $\f{(e_2)+e_4}{2}$ &
         $\frac{(e_2+e_3)+e_4}{2}$    \\& & & & & & &\\
$b_{\bf41}$     & $\frac{e_3}{2}$ & $\f{e_1+e_2}{2}$ &
       $\f{(e_1)+e_3+e_4}{2}$ & $\frac{(e_1)+e_2+(e_3)}{2}$ &
           $\f{e_1+e_4}{2}$ & $\f{(e_2+e_3)+e_4}{2}$  &
           $\frac{(e_2)+e_4}{2}$ \\ & & & & & & & \\
 $b_{\bf42}$     & $\frac{e_3+e_4}{2}$ & $\f{e_2+e_4}{2}$ &
        $\f{(e_1)+e_3}{2}$  & $\frac{e_2+(e_3)}{2}$ &
         $\f{e_1+e_4}{2}$ & $\f{(e_1+e_2+e_3)+e_4}{2}$ &
         $\frac{e_1+(e_2)}{2}$    \\ & & & & & & & \\
 \hline & & & & &  & & \\
  $\Lambda^B$ & $e_2, e_3, e_4$ & $e_1, e_2, e_4$ & $e_3,e_4$ & $e_2, e_4$ & $e_1, e_3, e_4$ & $e_4$ & $e_1,e_4$ \\ & & & & & & & \\
  $n_B$ & 3 & 3 & 2 & 2 & 3 & 1 & 2 \\ & & & & & & & \\
  $vol(\Lambda^B)$ & 1 & 1 & 1 & 1 & 1 & 1 &1
    \\ & & & & & & &\\ \hline & &  &  & & & & \\
  $tr_p(B)$ & $K_p^4(1)$  & $K_p^4(1)$  & $K_p^4(2)$ & $K_p^4(2)$  & $K_p^4(1)$  & $K_p^4(3)$ & $K_p^4(2)$
\\ & & & & & & &  \\  \hline
 \end{tabular}                          }
\end{center}

\medskip
\noindent
{\sl Table 2.18}: ${\mathcal F}_6  =\{{\bf 43}, {\bf 44}\}\, (6/2/1)$  ($F\simeq  \Z_2^3$, {\em diagonal, non-orientable}, $\beta_1=0, \beta_2=0$)

\begin{center}          \scriptsize{
 \begin{tabular}{|c||c|c|c|c|c|c|c|} \hline  & & & & & & &\\& $\g_1$ & $\g_2$ & $\g_3$ & $\g_1\g_2$ & $\g_1 \g_3$ & $\g_2 \g_3$ & $\g_1 \g_2 \g_3$ \\& & & & & & & \\  \hline & & & & & & & \\
  $B$  &  $d(I,-\tilde I)$  &  $d(\tilde I,-I)$  &  $d(-I,I)$ &  $d(\tilde I, \tilde I)$  &  $d(-I,-\tilde I)$  &  $d(-\tilde I,-I)$ & $d(-\tilde I,\tilde I)$ \\    & & & & & & & \\ \hline & & & & & & & \\
 $b_{\bf43}$     & $\frac{e_3}{2}$ & $\f{e_2}{2}$ &
           $\f{(e_1)+e_4}{2}$     & $\frac{e_2+(e_3)}{2}$ &
           $\f{(e_1)+e_3+(e_4)}{2}$ & $\f{e_1+(e_2+e_4)}{2}$  &
           $\frac{e_1+(e_2+e_3)+e_4}{2}$      \\ & & & & & & & \\
 $b_{\bf44}$        & $\frac{e_2+e_3}{2}$ & $\f{e_2+(e_4)}{2}$ &
        $\f{(e_1)+e_4}{2}$     & $\frac{(e_3)+e_4}{2}$ &
         $\f{(e_1+e_2)+e_3+(e_4)}{2}$ & $\f{e_1+(e_2)}{2}$ &
         $\frac{e_1+(e_3)}{2}$   \\ & & & & & & & \\
 \hline & & & & &  & & \\
  $\Lambda^B$ & $e_1, e_2, e_3$ & $e_2$ & $e_3,e_4$ & $e_2, e_4$ & $e_3$ & $e_1$ & $e_1,e_4$ \\ & & & & & & & \\
  $n_B$ & 3 & 1 & 2 & 2 & 1 & 1 & 2 \\ & & & & & & & \\
  $vol(\Lambda^B)$ & 1 & 1 & 1 & 1 & 1 & 1 & 1
    \\ & & & & & & &\\ \hline & &  &  & & & & \\
  $tr_p(B)$ & $K_p^4(1)$  & $K_p^4(3)$  & $K_p^4(2)$ & $K_p^4(2)$  & $K_p^4(3)$  & $K_p^4(3)$ & $K_p^4(2)$
\\ & & & & & & &  \\  \hline                  \end{tabular}             }
\end{center}

\medskip
\noindent
{\sl Table 2.19}: ${\mathcal F}_7={\bf 57, 58}\, (13/1/1)$  ($F\simeq  \Z_2 \times \Z_4$, {\em non-diagonal, non-orientable}, $\beta_1=1, \beta_2=1$)

\begin{center}    \scriptsize{
 \begin{tabular}{|c||c|c|c|c|c|c|c|} \hline  & & & & & & &\\& $\g_1$ & $\g_1^2$ & $\g_1^3$ & $\g_2$ & $\g_1 \g_2$ & $\g_1^2 \g_2$ & $\g_1^3 \g_2$ \\& & & & & & & \\  \hline & & & & & & & \\
  $B$  &  $d(\tilde I, -\tilde J)$  &  $d(I,-I)$  &  $d(\tilde I, \tilde J)$ &  $d(\tilde I, I)$  &  $d(I,-\tilde J)$  &  $d(\tilde I, -I)$ & $d(I,\tilde J)$ \\    & & & & & & & \\ \hline & & & & & & & \\
 $b_{\bf 57}$     & $\frac{e_2}{4}$ & $\f{e_2}{2}$ &
           $\f{3e_2}{4}$     & $\frac{e_3+e_4}{2}$ &
           $\f{e_2}{4}+\left( \frac{e_3+e_4}2\right)$ & $\f{e_2+(e_3+e_4)}{2}$  &
           $\frac{3e_2}{4}+ \left(\f{e_3 +e_4}2\right)$     \\ & & & & & & & \\
$b_{\bf 58}$     & $\frac{e_2}{4}$ & $\f{e_2}{2}$ &
           $\f{3e_2}{4}$     & $\!\frac{(-e_1)+e_3+e_4}{2}\!$ &
           $\!\f{e_2}{4}+\left(\frac{-e_1+e_3+e_4}2\right)\!$ & $\!\f{e_2+(-e_1+e_3+e_4)}{2}\!$  &
           $\!\frac{3e_2}{4}\!+\!\left( \f{e_1 - e_3 - e_4}2\right)\!$     \\ & & & & & & & \\     \hline & & & & &  & & \\
  $\Lambda^B$ & $e_2$ & $e_1, e_2$ & $e_2$ & $e_2, e_3, e_4$ & $e_1, e_2$ & $e_2$ & $e_1,e_2$ \\ & & & & & & & \\
  $n_B$ & 1 & 2 & 1 & 3 & 2 & 1 & 2 \\ & & & & & & & \\
  $vol(\Lambda^B)$ & 1 & 1 & 1  & 1  & 1  & 1  & 1  \\ & & & & & & &\\ \hline & &  &  & & & & \\
  $tr_p(B)$ & 1 0 0 0 -1  & $K_p^4(2)$  & 1 0 0 0 -1 & $K_p^4(1)$  & 1 2 2 2 1   & $K_p^4(3)$ & 1 2 2 2 1
\\ & & & & & & &  \\  \hline
 \end{tabular}                  }
\end{center}

\pagebreak
\medskip
\noindent
{\sl Table 2.20}: ${\bf 54}\, (12/3/4)$  ($F\simeq  D_4$, {\em non-diagonal, non-orientable}, $\beta_1=0, \beta_2=0$)
\

\begin{center}    \scriptsize{
 \begin{tabular}{|c||c|c|c|c|c|c|c|} \hline  & & & & & & &\\& $\g_1$ & $\g_1^2$ & $\g_1^3$ & $\g_2$ & $\g_1 \g_2$ & $\g_1^2 \g_2$ & $\g_1^3 \g_2$ \\& & & & & & & \\  \hline & & & & & & & \\
  $B$  &  $d(J,\tilde J)$  &  $d(I,-I)$  &  $d(J,-\tilde J)$ &  $d(J,-1,1)$  &  $d(I,J)$  &  $d(J,1,-1)$ & $d(I,-J)$ \\    & & & & & & & \\  & & & & & & & \\
 $b_{\bf54}$     & $\frac{e_1+(e_4)}{2}$ & $\f{e_1+e_2+(e_3+e_4)}{2}$ &
           $\f{e_2+(e_3)}{2}$     & $\frac{e_4}{2}$ &
           $-\f{e_2}{2}$ & $\f{e_1+e_2+e_3}{2}$  &
           $\frac{e_1+(e_3+e_4)}{2}$      \\ & & & & & & & \\     \hline & & & & &  & & \\
  $\Lambda^B$ & $e_1+e_2$ & $e_1, e_2$ & $e_1+e_2$ & $e_1+e_2, e_4$ & $e_1, e_2,e_3+e_4$ & $e_1+e_2, e_3$ & $e_1,e_2,e_3-e_4$ \\ & & & & & & & \\
  $n_B$ & 1 & 2 & 1 & 2 & 3 & 2 & 3 \\ & & & & & & & \\
  $vol(\Lambda^B)$ & $\sqrt{2}$ & 1 & $\sqrt{2}$  & $\sqrt{2}$  & $\sqrt{2}$  & $\sqrt{2}$  & $\sqrt{2}$
    \\ & & & & & & &\\ \hline & &  &  & & & & \\
  $tr_p(B)$ & 1 0 0 0 -1  & $K_p^4(2)$  & 1 0 0 0 -1 & $K_p^4(2)$  & $K_p^4(1)$  & $K_p^4(2)$ & $K_p^4(1)$
\\ & & & & & & &  \\  \hline
 \end{tabular}                  }
\end{center}

\medskip
\noindent
{\sl Table 2.21}: ${\bf 56}\, (12/4/3)$  ($F\simeq  D_4$, {\em non-diagonal, non-orientable}, $\beta_1=0, \beta_2=0$)
\

\begin{center}           \scriptsize{
 \begin{tabular}{|c||c|c|c|c|c|c|c|} \hline  & & & & & & &\\& $\g_1$ & $\g_1^2$ & $\g_1^3$ & $\g_2$ & $\g_1 \g_2$ & $\g_1^2 \g_2$ & $\g_1^3 \g_2$ \\& & & & & & & \\  \hline & & & & & & & \\
  $B$  &  $d(\tilde I,\tilde J)$  &  $d(I,-I)$  &  $d(\tilde I,-\tilde J)$ &  $d(-\tilde I,-J)$  &  $d(-I,\tilde I)$  &  $d(-\tilde I,J)$ & $d(-I,-\tilde I)$ \\    & & & & & & & \\  & & & & & & & \\
 $b_{\bf56}$     & $\frac{e_2}{4}+(\f{e_3}{2})$ & $\f{e_2+(e_3+e_4)}{2}$ &
           $\f{3e_2}{4}+(\f{e_4}{2})$     & $\frac{e_1}{2}$ &            $\big(\f{2e_1+e_2}{4}\big)+\f{e_4}{2}$ & $\f{e_1+(e_2)+e_3+e_4}{2}$  &
           $\big(\frac{2e_1+3e_2}4\big)+\f{e_3}{2}$      \\ & & & & & & & \\     \hline & & & & &  & & \\
  $\Lambda^B$ & $e_2$ & $e_1, e_2$ & $e_2$ & $e_1, e_3-e_4$ & $e_4$ & $e_1, e_3+e_4$ & $e_3$ \\ & & & & & & & \\
  $n_B$ & 1 & 2 & 1 & 2 & 1 & 2 & 1 \\ & & & & & & & \\
  $vol(\Lambda^B)$ & 1 & 1 & 1 & $\sqrt{2}$  & 1  & $\sqrt{2}$  & 1
    \\ & & & & & & &\\ \hline & &  &  & & & & \\
  $tr_p(B)$ & 1 0 0 0 -1  & $K_p^4(2)$  & 1 0 0 0 -1 & $K_p^4(2)$  & $K_p^4(3)$  & $K_p^4(2)$ & $K_p^4(3)$
\\ & & & & & & &  \\  \hline
 \end{tabular} }
\end{center}

%\pagebreak
%\medskip
\noindent
{\sl Table 2.22}: ${\mathcal F}_8=\{{\bf 60},{\bf 61},{\bf 62}\}\, (13/4/1)$  ($F\simeq  D_4$, {\em non-diagonal, orientable}, $\beta_1=1, \beta_2=0$)
\

\begin{center}               \scriptsize{
 \begin{tabular}{|c||c|c|c|c|c|c|c|} \hline  & & & & & & &\\& $\g_1$ & $\g_1^2$ & $\g_1^3$ & $\g_2$ & $\g_1 \g_2$ & $\g_1^2 \g_2$ & $\g_1^3 \g_2$ \\& & & & & & & \\  \hline & & & & & & & \\
  $B$  &  $d(I,\tilde J)$  &  $d(I,-I)$  &  $d(I,-\tilde J)$ &  $d(\tilde I,J)$  &  $d(\tilde I,-\tilde I)$  &  $d(\tilde I,-J)$ & $d(\tilde I,\tilde I)$ \\    & & & & & & & \\ \hline & & & & & & & \\
 $b_{\bf60}$     & $\frac{e_1}{4}$ & $\f{e_1}{2}$ &
           $\f{3e_1}{4}$     & $\frac{e_2}{2}$ &            $(\!\f{e_1}{4}\!)\!+\!\f{e_2}{2}$ & $\f{(\!e_1\!)\!+\!e_2}{2}$  &
           $(\!\frac{3e_1}{4}\!)\!+\!\f{e_2}{2}$      \\ & & & & & & & \\
 $b_{\bf61}$     & $\frac{e_1}{4}\!+\!(\f{e_4}{2})$ & $\f{e_1+(e_3+e_4)}{2}$ &       $\f{3e_1}{4}\!-\!(\f{e_3}{2})$     & $\frac{e_2}{2}$ &      $(\!\f{e_1}{4}\!)\!+\!\f{e_2+e_3}{2}$ &
$\f{(\!e_1\!)\!+\!e_2\!+\!e_3\!+\!e_4}{2}$  &
           $(\!\frac{3e_1}{4}\!)\!+\!\f{e_2+e_4}{2}$      \\ & & & & & & & \\
 $b_{\bf62}$     & $\frac{e_1}{4}\!+\!\f{e_2+(e_4)}{2}$ & $\f{e_1+(e_3+e_4)}{2}$ &
           $\f{3e_1}{4}\!+\!\f{e_2+(e_3)}{2}$     & $\frac{e_2}{2}$ &            $(\!\f{e_1}{4}\!)\!+\!\f{e_3}{2}$ & $\f{(\!e_1\!)\!+\!e_2\!+\!e_3\!+\!e_4}{2}$  &
           $(\!\frac{3e_1}{4}\!)\!+\!\f{e_4}{2}$      \\ & & & & & & & \\
\hline & & & & &  & & \\
  $\Lambda^B$ & $e_1,e_2$ & $e_1, e_2$ & $e_1,e_2$ & $e_2, e_3\!+\!e_4$ & $e_2,e_3$ & $e_2, e_3\!-\!e_4$ & $e_2,e_4$ \\ & & & & & & & \\
  $n_B$ & 2 & 2 & 2 & 2 & 2 & 2 & 2 \\ & & & & & & & \\
  $vol(\Lambda^B)$ & 1 & 1 & 1 & $\sqrt{2}$  & 1  & $\sqrt{2}$  & 1
    \\ & & & & & & &\\ \hline & &  &  & & & & \\
  $tr_p(B)$ & 1 2 2 2 1  & $K_p^4(2)$  & 1 2 2 2 1 & $K_p^4(2)$  & $K_p^4(2)$  & $K_p^4(2)$ & $K_p^4(2)$
\\ & & & & & & &  \\  \hline
 \end{tabular}                      }
\end{center}

%\subsubsection*{Example: A $1$-isospectral pair with holonomy group $D_4$.}
\begin{ejem}[{\it A $1$-isospectral pair with holonomy group $D_4$}]
As an application of the methods in this section, we will use formula
(\ref{multip}) to carry out explicit calculations of multiplicities of
eigenvalues for the groups ${\bf 60}$ and ${\bf 61}$, both having holonomy
groups isomorphic to $D_4$. As a consequence, we shall see that ${\bf
60,61}$ are 1-isospectral and 3-isospectral but they are not
$p$-isospectral for $p\not=1,3$. Since the manifolds associated to the
groups ${\bf 60}, {\bf 61}$ are orientable, 1-isospectral and
3-isospectral are equivalent in this case.

According to formula~(\ref{multip}) we have:
\begin{equation}            \label{e.dp6061}
d_{p,\mu}(\G)=\tf18\tbinom{4}{p}|\Ld_\mu| + \tf18 \sum_{\begin{smallmatrix} \gamma=BL_b \in \Ld \backslash \G \\ \g \not=\I \end{smallmatrix}} tr_p(B)\; e_{\mu,\g}(\G)
\end{equation}
where $\Ld_\mu=\{\ld\in\Ld : ||\ld||^2=\mu \}$.
For the elements $\g= \g_1^2,\g_2,\g_1\g_2,\g_1^2\g_2,\g_1^3\g_2$,  we have that $tr_p(B)=K_p^4(2)$ (see Table~2.22 above), so for $p=1$ or $p=3$ these traces vanish (see (\ref{krawtvalues})).
Thus, using that $tr_p(B_1)=tr_p(B_1^3)=2$ for $p=1,3$, we get
\begin{equation}                \label{e.d16061}
d_{1,\mu}(\G)=d_{3,\mu}(\G) =\tf12|\Ld_\mu| + \tf14 e_{\mu, \g_1}(\G) +
\tf14 e_{\mu, \g_1^3}(\G).
\end{equation}
In order to check 1-isospectrality it suffices to show that
$e_{\mu, \g_1}(\G_{\bf 60})=e_{\mu, \g_1}(\G_{\bf 61})$ and
$e_{\mu, \g_1^3}(\G_{\bf 60})=e_{\mu, \g_1^3}(\G_{\bf 61})$. We
have:
$$ e_{\mu,\g_1}(\G_{\bf 60})=\sum_{ \begin{smallmatrix} v\in \Ld_\mu \\ B_1 v=v \end{smallmatrix} } e^{ -2\pi i  v \cdot \f{e_1}{4}} = \sum_{ \begin{smallmatrix} v\in \Z e_1 \oplus \Z e_2 \\ ||v||^2=\mu \end{smallmatrix}} e^{-\f{\pi }{2} i   v\cdot e_1}$$
$$ e_{\mu,\g_1}(\G_{\bf 61})=\sum_{ \begin{smallmatrix} v\in \Ld_\mu \\ B_1 v=v \end{smallmatrix} } e^{ -2\pi i  v\cdot (\f{e_1}{4}+\f{e_4}{2})} = \sum_{ \begin{smallmatrix} v\in \Z e_1 \oplus \Z e_2 \\ ||v||^2=\mu \end{smallmatrix}} e^{-\f{\pi }{2} i   v \cdot e_1 } e^{-\pi i  v\cdot e_4 }.$$
Since $e^{-\pi i  v \cdot e_4}=1$ for $v\in \Z e_1\oplus \Z e_2$,  we have  that $e_{\mu, \g_1}(\G_{\bf 60})=e_{\mu, \g_1}(\G_{\bf 61})$. Similarly, we get that  $e_{\mu, \g_1^3}(\G_{\bf 60})=e_{\mu, \g_1^3}(\G_{\bf 61})$.
Observe that we can read off this fact from the rows of the vectors $b$ in Table 2.22, since we know that the orthogonal projection of the vectors in between parentheses onto the space fixed by the corresponding matrix $B$, is zero. In this way we have proved that the manifolds $M_{\bf 60}, M_{\bf 61}$ obtained from the groups ${\bf 60, 61}$ are 1-isospectral and hence 3-isospectral.

We now show that they can not be $0,2$ nor 4-isospectral. Take $\mu=1$. The formula for the multiplicities for $p=0,2,4$ in this particular case is:
$$ d_{p,1} = \tf18 \tbinom{4}{p}|\Ld_1| + \tf18 tr_p(B_1)\, (e_{1, \g_1} + e_{1, \g_1^3}) +
\tf18 K_4^p(2)\,(e_{1, \g_1^2} + e_{1, \g_2} + e_{1, \g_1\g_2} + e_{1, \g_1^2\g_2} + e_{1, \g_1^3\g_2} ).$$

Computing  the corresponding $e_{1,\g}$'s we obtain

\begin{center}
\begin{tabular}{|c|c|c|c|c|c|c|c|} \hline &&&&&&& \\
$\g$ & $\g_1$ & $\g_1^2$ & $\g_1^3$ & $\g_2$ & $\g_1\g_2$ & $\g_1^2\g_2$ & $\g_1^3\g_2$ \\ &&&&&&& \\ \hline  & & & & & & & \\
$e_{1, \g}(\G_{\bf 60})$ & 2& 0& 2& -2& 0& -2& 0 \\ &&&&&&& \\ \hline &&&&&&& \\
$e_{1, \g}(\G_{\bf 61})$ & 2& 0& 2& -2& -4& -2& -4\\ &&&&&&& \\ \hline
\end{tabular}
\end{center}

\smallskip
>From this, and taking into acount that $\Ld_1=\{\pm e_1; \pm e_2; \pm e_3; \pm e_4\}$,  we get
\begin{equation*}
\begin{split}
 d_{p,1}(\G_{\bf 60})= \tbinom{4}{p} + \tf12 tr_p(B_1) - \tf12 K_p^4(2) \\
 d_{p,1}(\G_{\bf 61})= \tbinom{4}{p} + \tf12 tr_p(B_1) - \tf32 K_p^4(2).
\end{split}
\end{equation*}

Since $K_p^4(2)\not=0$ for $p=0,2,4$ (see (\ref{krawtvalues})), it is clear that $  d_{p,1}(\G_{\bf 60})\not= d_{p,1}(\G_{\bf 61})$
and this shows that $M_{\bf 60}$ and $M_{\bf 61}$ can not be $p$-isospectral for $p=0,2,4$.
   For example, for $p=0,4$, we have $d_{p,1}(\G_{\bf 60})=1$ and  $d_{p,1}(\G_{\bf 61})=0$, while for $p=2$ we have
$d_{2,1}(\G_{\bf 60})=6$ and  $d_{2,1}(\G_{\bf 61})=4$.

\end{ejem}

\bigskip

\subsubsection*{Appendix.}
Here we list the matrices $X\ne\I$ and $C$, corresponding to the
groups described in the tables in this section.
%\medskip
$$ X_{\bf 3}= X_{\bf 12}=\left ( \begin{smallmatrix}
1 & 0 & 0 & 0\\     0 & 1 & 0 & 0\\ 0 & 0 & 1 & 0\\ 0 & 0 & 1 & 1
\end{smallmatrix} \right ),\; X_{\bf 6}=\left (\begin{smallmatrix}
0 & 1 & 0 & 0\\     1 & 0 & 0 & 1\\ 1 & 0 & 0 & 0\\ 0 & 0 & 1 & 0
\end{smallmatrix} \right ),\; X_{\bf 13}=X_{\bf 14}=X_{\bf 15}
=X_{\bf 22} =\left (\begin{smallmatrix} 0 & 1 & 0 & 0\\     1 & 0
& 0 & 0\\ 1 & 0 & 1 & 0\\ 0 & 0 & 0 & 1   \end{smallmatrix} \right
),$$
$$X_{\bf 28}=X_{\bf 29}=\left (\begin{smallmatrix}
0 & 0 & 0 & 1\\     1 & 0 & 0 & 0\\ 0 & -1 & 0 & 0\\ 0 & 0 & 1 & 0
\end{smallmatrix} \right ), \; X_{\bf 47}=\left
(\begin{smallmatrix} 1 & 0 & 0 & 0\\     0 & -1 & 1 & 0\\ 0 & -1 &
0 & -1\\ 0 & -1 & 0 & 0   \end{smallmatrix} \right ).$$

$$ C_{\bf 2}'= C_{\bf 29}'=\left ( \begin{smallmatrix}
1 & & & \\ & 1 & 1 & \\ & & 1 & \\ & & & 1 \end{smallmatrix}
\right ),
 C_{\bf 2}''= \left ( \begin{smallmatrix}
1 & & 1 & \\     & 1 & 1 & \\  & & 1 & \\ & & & 1
\end{smallmatrix} \right ),
 C_{\bf 3}'= C_{\bf 9}'= C_{\bf 10}'=C_{\bf 11}'=\left ( \begin{smallmatrix}
1 & & & \\  1 & 1 &  & \\  & & 1 & \\ & & & 1 \end{smallmatrix}
\right ),
 C_{\bf 3}''= \left ( \begin{smallmatrix}
1 & & & \\   & 1 &  & \\ 1 & & 1 & \\ 1 & & & 1 \end{smallmatrix}
\right ), $$
$$ C_{\bf 3}'''=  C_{\bf 47}' = \left ( \begin{smallmatrix}
1 & & & \\  1 & 1 &  & \\ 1 & & 1 & \\ 1 & & & 1
\end{smallmatrix} \right ),
 C_{\bf 5}'= \left ( \begin{smallmatrix}
1 & & & \\   & 1 &  & \\  & & 1 & 1 \\ & & & 1  \end{smallmatrix}
\right ),
 C_{\bf 6}'= C_{\bf 13}'=C_{\bf 14}'=C_{\bf 15}'=C_{\bf 22}'=\left ( \begin{smallmatrix}
1 & & & \\  1 & 1 & & \\ 1 & & 1 & \\ & & & 1 \end{smallmatrix}
\right ), $$
\begin{tabular}{c} \qquad $ C_{\bf 7}'= C_{\bf 8}'=C_{\bf 10}''= C_{\bf 12}'=
C_{\bf 19}'= C_{\bf 20}'= C_{\bf 21}'= C_{\bf 23}'= C_{\bf 45}'=
\left ( \begin{smallmatrix}
        1 & 1 &   &   \\
         & 1 &   &   \\
          &   & 1 &   \\
          &   &   & 1 \\ \end{smallmatrix} \right ). $ \\ \\ \\ \\ \\ \qquad
\end{tabular}

%\bigskip

Finally, we shall give here a list with the Sunada numbers
$c_{d,t}$'s for those groups having diagonal holonomy
representation. We only give the values for $c_{1,1}, c_{2,1},
c_{2,2}, c_{3,1}, c_{3,2}, c_{3,3}$, since
$c_{1,0}=c_{2,0}=c_{3,0}=0$, $c_{4,0}=1$ and
$c_{4,1}=c_{4,2}=c_{4,3}=c_{4,4}=0$ for any 4-dimensional
Bieberbach group of diagonal type.
\begin{center}

\renewcommand{\arraystretch}{.25}
\noindent \begin{tabular}{|l|c|cc|ccc|}
    \hline &&&&&& \\ $\Gamma_\#$ &  $c_{1,1}$ & $c_{2,1}$ & $c_{2,2}$ &  $c_{3,1}$ & $c_{3,2}$ & $c_{3,3}$ %& $c_{4,0}$
\\ &&&&&& \\ \hline  &&&&&& \\
$\bf 2$   & 0 & 0 & 0 & 1 & 0 & 0 %& 1
\\
$\bf 2'$  & 0 & 0 & 0 & 0 & 1 & 0 %& 1
\\
$\bf 2''$ & 0 & 0 & 0 & 0 & 0 & 1 %& 1
\\
$\bf 4$   & 1 & 0 & 0 & 0 & 0 & 0 %& 1
\\
$\bf 5$   & 0 & 1 & 0 & 0 & 0 & 0 %& 1
\\
$\bf 5'$  & 0 & 0 & 1 & 0 & 0 & 0 %& 1
\\ &&&&&& \\ \hline
&&&&&& \\
$\bf 7$   & 0 & 1 & 0 & 2 & 0 & 0 %& 1
\\
$\bf 7'$  & 0 & 0 & 1 & 1 & 1 & 0 %& 1
\\
$\bf 8$   & 0 & 1 & 0 & 1 & 1 & 0 %& 1
\\
$\bf 8'$  & 0 & 0 & 1 & 1 & 0 & 1 %& 1
\\
$\bf 9$   & 0 & 0 & 1 & 2 & 0 & 0 %& 1
\\
$\bf 9'$  & 0 & 1 & 0 & 1 & 1 & 0 %& 1
\\
$\bf10$   & 0 & 0 & 1 & 1 & 1 & 0 %& 1
\\
$\bf10'$  & 0 & 1 & 0 & 1 & 0 & 1 %& 1
\\
$\bf10''$ & 0 & 1 & 0 & 0 & 2 & 0 %& 1
\\
$\bf11$   & 0 & 0 & 1 & 0 & 2 & 0 %& 1
\\
$\bf11'$  & 0 & 1 & 0 & 0 & 1 & 1 %& 1
\\ &&&&&& \\ \hline
&&&&&& \\
$\bf18$  & 1 & 1 & 0 & 1 & 0 & 0 %& 1
\\
$\bf19$  & 1 & 1 & 0 & 1 & 0 & 0 %& 1
\\
$\bf19'$ & 1 & 1 & 0 & 0 & 1 & 0 %& 1
\\
$\bf20$  & 1 & 0 & 1 & 1 & 0 & 0 %& 1
\\
$\bf20'$ & 1 & 0 & 1 & 0 & 1 & 0 %& 1
\\
$\bf21$  & 1 & 1 & 0 & 0 & 1 & 0 %& 1
\\
$\bf21'$ & 1 & 1 & 0 & 0 & 0 & 1 %& 1
\\
&&&&&& \\ \hline
\end{tabular}
  \renewcommand{\arraystretch}{0.53}
\begin{tabular}{|l||c|cc|ccc|}
    \hline &&&&&& \\ $\Gamma_\#$ &  $c_{1,1}$ & $c_{2,1}$ & $c_{2,2}$ &  $c_{3,1}$ & $c_{3,2}$ & $c_{3,3}$ %& $c_{4,0}$
\\ &&&&&& \\ \hline  &&&&&& \\
$\bf23$ & 2 & 1 & 0 & 0 & 0 & 0 \\
$\bf23'$ & 2 & 0 & 1 & 0 & 0 & 0  \\
 &&&&&& \\ \hline  &&&&&& \\
$\bf24$ & 0 & 3 & 0 & 0 & 0 & 0  \\
$\bf25$ & 0 & 2 & 1 & 0 & 0 & 0  \\
$\bf26$ & 0 & 3 & 0 & 0 & 0 & 0  \\
$\bf27$ & 0 & 1 & 2 & 0 & 0 & 0  \\ &&&&&& \\ \hline
&&&&&& \\
$\bf33$ & 1 & 2 & 1 & 3 & 0 & 0  \\
$\bf34$ & 1 & 1 & 2 & 2 & 1 & 0  \\
$\bf35$ & 1 & 3 & 0 & 1 & 1 & 1  \\
$\bf36$ & 1 & 3 & 0 & 2 & 1 & 0  \\
$\bf37$ & 1 & 2 & 1 & 2 & 0 & 1  \\
$\bf38$ & 1 & 1 & 2 & 2 & 1 & 0 \\
$\bf39$ & 1 & 2 & 1 & 1 & 2 & 0  \\
$\bf40$ & 1 & 3 & 3 & 1 & 1 & 1  \\
$\bf41$ & 1 & 2 & 1 & 1 & 2 & 0  \\
$\bf42$ & 1 & 3 & 0 & 0 & 3 & 0  \\ &&&&&& \\ \hline
&&&&&& \\
$\bf43$ & 3 & 2 & 1 & 1 & 0 & 0  \\
$\bf44$ & 3 & 3 & 0 & 0 & 1 & 0  \\ &&&&&&\\ \hline
 \end{tabular}
\end{center}

\section{Zeta functions and $p$-heat trace polynomials} \label{s.sec3}
    In this section we will show that the $p$-heat trace $Z_p^{\G}(s)$
    of an $n$-dimensional  compact flat manifold $M_\G$ can be viewed
    as a polynomial of degree $n$ in
    a finite number of algebraically independent functions
    $x(s),y(s),z_1(s), \dots, z_m(s)$,  for $s>0$.
We let the dimension $n$ of $M_\G$ be arbitrary, until further
notice.

\sk We begin by  recalling some facts on the $L$-spectrum of a
Bieberbach manifold $\G$. For $\g=BL_b$, denote by $p_B : \R^n
\rightarrow (\R^n)^B$ the orthogonal projection and put
$b_+=p_B(b)$. The $L$-spectrum of $M_\G$  is given by
$Spec_{L}(M_\G)=\{||b_+ + \ld_+|| : BL_b \in \G \}$ with $BL_b$
running through a full set of representatives of $F$ and $\ld \in
\Ld$ (see Section 1).
 Note that for the $L$-spectrum, only the components of $b$ not between parentheses in the tables count. %Let look at the group numbered ${\bf 2}$ in Table 1. In that case we have $\Ld^B=\langle e_1,e_2,e_3 \rangle$ and $b=\f{1}{2}e_1$ so $b_+=b$.
In this section, to each element $B$ we will associate a monomial in several variables and will express the corresponding zeta functions in terms of these monomials, thus getting polynomial expressions for the $p$-heat traces in which the $L$-spectrum information is encoded. We will give a complete list of these polynomials from which one can recover the $L$-spectrum of the corresponding manifold. % rather easily.

By Theorem \ref{t.zetafunc}, for a compact flat manifold $M_\Gamma = \vcp$ with translation lattice $\Lambda$ and holonomy group $F$ we have the following expression for the $p$-heat traces:
\begin{equation}                            \label{eq.longzetaf2}
         Z^{\,\Gamma}_p(s) = \tfrac{1}{|F|} \sum_{BL_b \in \Ld\backslash\G}         \frac{\text{tr}_p(B)}{\text{vol}({\Lambda ^*}^{\text{B}})}\,
        (4 \pi s)^{-\f{n_B}2}
        \sum _{\lambda_+ \in p_B(\Ld)}
        e^{- \frac{ {\|b_+ + \lambda_+ \|}^2}{4s}}.
    \end{equation}

and if $\Gamma$ is of diagonal type, with $F\simeq \Z_2^r$, we can write:
\begin{equation}                            \label{eq.diagzetaf2}
     Z^{\,\Gamma}_p(s)=
        \frac{1}{2^r} \sum _{d=1}^n K_p^n(n-d)\,
        (4 \pi s)^{-\frac{d}2} \sum _{t=0}^d c_{d,t}(\G) \: \theta_{d,t}
        (\tfrac 1{4s})
    \end{equation}
where, for $Re(s)>0$,
    \begin{equation*}                           \label{eq.tethas2}          \theta_{d,t}(s) = \sum
        _{(m_1,\dots,m_d)\in\Z^d}%{\Sb m_j \in \Z\\ 1\le j \le d \endSb}
        e^{-s \left(\sum_{j=1}^t (\frac 12 + m_j)^2 + \sum_{j=t+1}^d        {m_j}^2\right)}.
    \end{equation*}

The theta functions $\theta_{d,t}$ have simple expressions in terms of $\theta_0$ and
$\theta_1$, where
\begin{equation}\label{thetas}
    \theta_0(s):=\theta_{1,0} (s) = \sum
        _{m\in\Z}e^{-s m^2},\qquad
 \theta_1 (s):=\theta_{1,1}(s)= \sum _{m\in\Z}e^{-s (\frac 12 + m)^2 }, \quad \text{ for } Re (s) >0.
    \end{equation}
We define the functions
$$x(s):= \frac{\theta_0(\frac1{4s})}{\sqrt{4\pi s}}
\qquad  \text{ and }  \qquad y(s):= \frac{\theta_1(\frac1{4s})}{\sqrt{4\pi s}}.$$
\begin{lema}\label{algindep} For $Re(s)>0$ we have $\theta_{d,t}(s) = \theta_0^{d-t}(s)
 \,\theta_1^t(s)$.
Furthermore, if $a_{d,t} \in \C$ and $\sum _{0\le t \le d} a_{d,t}
x(s)^{d-t}y(s)^t =0$, then $a_{d,t}=0$ for any $0\le d\le t$.
\end{lema}
\begin{proof}
Decomposing $\Z^d$ as $\Z^{d-t} \oplus \Z^t$ we get
    $$\theta_{d,t}(s) = \sum_{\scriptsize (m_1,\dots, m_t) \in \Z^t}
        e^{-s \sum_{j=1}^t (\frac 12 + m_j)^2} \;
        \sum_{\scriptsize (m_1,\dots,m_{d-t}) \in \Z^{d-t}}
        e^{-s  \sum_{j=1}^{d-t}  {m_j}^2 } = \theta_0^{d-t}(s) \theta_1^t(s).$$

To prove the second assertion, we note  that, as $s\rightarrow
+\infty$,
$$ \theta_0^{t}(s)\sim 1, \qquad  \theta_1^t(s)\sim 2^t e^{-\frac {ts}
4}$$ Now, if $u=1/4s$ then, as $s \downarrow 0$,
\begin{eqnarray*}
\sum _{0\le t \le d} a_{d,t}\, x(s)^{d-t}y(s)^t &=& \sum _{0\le t
\le d} a_{d,t}\;(\tfrac u \pi)^{d/2}\theta_0(u)^{d-t}
\theta_1(u)^t
\\& \sim & \sum _{0\le t \le d} a_{d,t}\;2^t(\tfrac u \pi)^{d/2}
e^{-\frac{tu} 4}
\end{eqnarray*}

Now if not all $a_{d,t}=0$, let $t_m$ be minimal with the property
that some $a_{d,t_m}\ne 0$ and let $d_M$ be  maximal among the $d$
so that $a_{d,t_m}\ne 0$. Multiplying the previous expression by
$u^{-d_M/2} e^{\frac{t_m u} 4}$ and letting $u\rightarrow \infty$
we have that all terms, except the one corresponding to $d=d_M,
t=t_m$, tend to zero. Thus, we get that $a_{d_M,t_m}2^t
\pi^{d/2}=0$, a contradiction.
\end{proof}

%We have the following:
\begin{prop}            \label{p.zpoly}
Let $\G$ be an $n$-dimensional Bieberbach group of diagonal type with canonical
 translation lattice $\Ld$ and holonomy group $F$.
 Then $Z^{\,\Gamma}_p(s) \in \frac1{|F|} \Z[x(s), y(s)]$  has degree $n$, no independent term,
  and the coefficient of $x(s)^k$ is 0 (resp.\@ 1) for all $k=1, \dots, n-1$
 (resp. $k=n$).
\end{prop}

\begin{proof}
By (\ref{eq.diagzetaf2}) and using the notation $Z_{\I} (s) := (4\pi s)^{-\tfrac n2} \sum_{\ld \in \Ld}  e^{-\frac{||\ld||^2}{4s}}$    we can write
 $$  Z^{\,\Gamma}_p(s) = \tfrac{\binom{n}{p}}{|F|} \, Z_{\I}(s) +
    \tfrac1{|F|}    \sum_{d=1}^{n-1}   K_p^n(n-d)   \sum_{t=0}^d
    c_{d,t}(\G) \: \frac{\theta_{d,t}(\f1{4s})}{(4\pi s)^{\f{d}{2}}
} .$$
But
$$Z_{\I} (s) = {(4\pi s)^{-\tfrac{n}{2}}} \sum_{(m_1,\dots, m_n) \in \Z^n}
e^{-\frac{\big(m_1^2+\cdots+m_n^2\big)}{4s}} =
\f{\theta_{n,0}(\tfrac{1}{4s})}{(\sqrt{4\pi s})^n}  = x(s)^n $$
and similarly,
\begin{equation}                            \label{eq.zetapoly}
   \f{\theta_{d,t}(\tfrac{1}{4s})}{({4\pi s})^{d/2}} =
    \f{\theta_{0}^{d-t}(\tfrac{1}{4s})}{(\sqrt{4\pi s})^{d-t}}      \f{\theta_{1}^t(\tfrac{1}{4s})}{(\sqrt{4\pi s})^t} = x(s)^{d-t}\,y(s)^t
\end{equation}
thus $$|F|\; Z^{\,\Gamma}_p(s) = \tbinom{n}{p} \, x(s)^n +
    \sum_{d=1}^{n-1}   K_p^n(n-d)   \sum_{t=0}^d    c_{d,t}(\G)\; x(s)^{d-t}\,y(s)^t .$$
Since all  numbers $\binom np$, $K_p^n(n-d)$ and $c_{d,t}(\G)$ are integers, we are done.
\end{proof}

In the case of groups with holonomy representation that is not of
diagonal type, one has an analogue of (\ref{eq.diagzetaf2}) but
the expression is more complicated. Since $\Ld$ is the canonical
lattice, $\Ld=\Ld^*$. Furthermore, we have that $1 \leq n_B \leq
n$ and $n_B=n$ if and only if $B=\I$. Thus, adding over $\g=BL_b$
with fixed $n_B=d$, we may write:
\begin{equation}                            \label{eq.varzetaf}
   Z^{\,\Gamma}_p(s) = \tfrac{\binom{n}{p}}{|F|} \,
    Z_{\I}(s) +     \tfrac1{|F|}
    \sum_{d=1}^{n-1} \,(4\pi s)^{-\f{d}{2}}
    \sum_{\scriptsize \begin{array}{c} \g= BL_b \in \Ld\backslash\G \\ n_B=d                         \end{array}}
    \frac{\text{tr}_p(B)}{\text{vol}({\Lambda}^{\text{B}})} \: Z_\g^\G (s) \end{equation}
where  $Z_\g^\G(s):=\sum _{\lambda_+ \in p_B(\Ld)}  e^{- \frac{
{\| \lambda_+ + b_+\|}^2}{4s}}$.

\begin{rem}             \label{p.zpoly2}
Let $\G$ a Bieberbach group of even (resp.\@ odd) dimension $n$
with diagonal holonomy representation and canonical lattice of
translations $\Ld$. If $M_{\G}$ is the associated compact flat
manifold then
 $M_\G$ is orientable if an only if  $Z_p^{\G}(x(s),y(s))$
has only even (resp.\@ odd) degree terms for all $0\le p\le n$, that is,
$c_{d,t}(\G)=0$ for all $d=1, 3, \dots, n-1$ (resp.\@ $d=0, 2, \dots, n$)
and for all $t$.

Indeed, let $r:\G \arr \text{O}(n)$ be the projection $r(BL_b)=B$.
It is known that $M_\G$ is orientable if and only if $r(\g) \in
\text{SO}(n)$ for all $\g \in \G$, that is, $n_B \in 2\Z$ (resp.\@
$n_B \in 2\Z +1$) for all $B\in F$ since $B$ is a diagonal matrix
for each $BL_b \in \G$. Now, using the expression
(\ref{eq.varzetaf}) with canonical $\Ld$ and diagonal matrices,
and comparing it with the expression (\ref{eq.diagzetaf2}), we see
that:
$$ \sum_{\scriptsize \begin{array}{c} \g = BL_b \in \Ld\backslash\G \\ n_B=d                                \end{array}} \text{tr}_p(B) (4\pi s)^{-\f{d}{2}} \;  Z_\g^\G(s)
         = \sum_{1\leq t \leq d}  K_p^n(n-d)    c_{d,t}(\G)\; x(s)^{d-t}\,y(s)^t  $$
are the terms of degree $d$, and this says  that $M_\G$ is orientable if and only if
 $Z_p^\G(x(s),y(s))$ has only even (resp.\@ odd) degree terms.
 Now since $$ \text{Coef}(x(s)^{d-t}y(s)^t)=K_p^n(n-d)\,c_{d,t}(\G)$$
 and  $K_p^n(n-d)\not= 0$  in general, it follows that $M_\G$ is orientable if and only if
 $c_{d,t}(\G)=0$ for each $d=1, 3, \dots, n-1$ if $n$ is even (resp.\@  $d=2, 4, \dots, n-1$
 if $n$ is odd), and for all $t$.
\end{rem}

\subsubsection*{Dimension 4}
We now look at the special case of dimension 4. From
Section 2 we know that the information of the $L$-spectrum is
encoded in the terms $Z_\g^\G (s)$ of $Z_p^\G(s)$. From expression
(\ref{eq.diagzetaf2}) we just have to look at the
$\theta_{d,t}$'s, or as we have seen, at the monomials
$x(s)^{d-t}\,y(s)^t$.
 In general, however, if $\G$ is not of diagonal type, the squared lengths  of the closed
 geodesics,  $\|\ld_+ + b_+\|^2$,
have summands different from $m^2$ or $(m+\f12)^2$. For $Re(s)>0$
and $r\in \frac1{|F|}\Z$, $0\le r<1$, we define the function
$$ \phi_{d,r}(s) := \sum_{(m_1, \dots, m_d) \in \Z^d }
e^{-s \big(\frac{m_1 + \cdots + m_d +r}{\sqrt d}\big)^2}$$ and we set
$$z_{d,r}(s):= \frac{ \phi_{d, r}(\tf{1}{4s})}{\sqrt{4\pi s}}.$$
\begin{rem}
Note that $\phi_{1,0}(z)=\theta_{0}(z)$ and
$\phi_{1,\scriptsize{\f{1}{2}}}(z)=\theta_{1}(z)$ so we could just
use the $\phi_{d,r}$'s, instead of the $\theta_{d,t}$'s, to
express the zeta functions as polynomials.
\end{rem}

\begin{prop}                        \label{p.zpoly3}
Let $\G$ be a 4-dimensional Bieberbach group with canonical
lattice  $\Ld$ and holonomy group $F$. Let $K=\Z(\sqrt 2, \sqrt
3)$. Then $Z_p^\G(s)$ is given by a polynomial expression with
coefficients in $K$ in the variables $x(s), y(s), z_{1,r}(s),
z_{2,r}(s), z_{3,r}(s)$, where $r=k/|F|$ and $0 \le k \le |F|-1$.
\end{prop}

\begin{proof} The proposition follows by case by case verification, directly
from the tables in Section 2. We need only check those groups with
non-diagonal holonomy representation. We will carry out the
verification  for just one group to ilustrate the method.

Take the group numbered~ ${\bf 29}$ (see Table 2.11). For the
first element $B_1$ we have $b_1=-e_2+e_3+e_4$ and $\Ld^{B_1}$ has
the basis  $\{e_2+e_3,e_4\}$. Since the vector $-e_2+e_3$ is
orthogonal to $e_2+e_3$, it appears in the table in between
parentheses. Thus $({b_1})_+=\tf{1}{2}e_4$. Take
$\ld=m_1e_1+\cdots +m_4e_4$ in $\Ld$. The orthogonal projection of
$\ld$ onto the subspace generated by $e_2+e_3$ is given by:
$$  \left((m_1e_1+m_2e_2+m_3e_3+m_4e_4)\cdot \tf{(e_2+e_3)}{\sqrt{2}}\right)
\tf{e_2+e_3} {\sqrt{2}} =
\tf{(m_2+m_3)}{\sqrt{2}}\tf{(e_2+e_3)}{\sqrt{2}},$$ thus ${b_1}_+
+ \ld_+ = (\f{m_2+m_3}{\sqrt 2}) \f{(e_2+e_3)}{\sqrt 2} +
(m_4+\tf{1}{2}) e_4$. In this way, for $\g_1=B_1L_{b_1}$, we have:
\begin{equation*}
\begin{split}
Z_{\g_1}^{\bf{29}}(s) &  =  \sum_{\ld_+ \in \Ld^{B_1}} e^{-
\tf{||{b_1}_+ + \ld_+||^2}{4s}  }  =  \sum_{(m_1,m_2,m_3) \in
\Z^3} e^{-\left\{ \tf{(\f{m_1+m_2}{\sqrt 2})^2 + (m_3+ \f{1}{2})^2
}{4s} \right\}}
 \\  &  =  \sum_{(m_1,m_2) \in \Z^2} e^{- \tf{(\f{m_1+m_2}{\sqrt 2})^2}{4s} }
  \sum_{m_3 \in \Z} e^{-\tf{({m_3+\f 12})^2}{4s}}   =
\phi_{2,0}(\tf{1}{4s}) \; \theta_{1}(\tf{1}{4s}).
\end{split}
\end{equation*}

Doing the same for $\g_2=B_2L_{b_2}$ and $\g_3=\g_1\g_2$ we get
$Z_{\g_2}^{\bf{29}}(s) = \phi_{2,0}(\tf{1}{4s}) \;
\theta_{1}(\tf{1}{4s})$ and $Z_{\g_3}^{\bf{29}}(s) =
\theta_{2,2}(\tf{1}{4s}) = \theta_{1}^2(\tf{1}{4s})$. Now, using
formula (\ref{eq.varzetaf}) we get for this case
\begin{equation*}
\begin{split}
Z^{\,\bf{29}}_p(s) & = \tfrac{\binom{4}{p}}{4} \,
    Z_{\I}(s) +     \tfrac{K_p^4(2)}{2} \frac{1}{4\pi s}
\left( \tf{\sqrt 2}2 Z_{\g_1}^{\bf{29}}(s) + \tf{\sqrt 2}2
Z_{\g_2}^{\bf{29}}(s) +
Z_{\g_3}^{\bf{29}}(s) \right)  \\
 & = \tfrac{\binom{4}{p}}{4} \,
    \;x(s)^4 +  \tfrac{K_p^4(2)}{2} \frac{1}{4\pi s}
\left( {\sqrt 2} \, \phi_{2,0}(\tf{1}{4s}) \; \theta_1(\tf{1}{4s})
+ \theta_1^2(\tf{1}{4s}) \right)  \\
& =   \tfrac{\binom{4}{p}}{4} \, \;x(s)^4 +
\tfrac{K_p^4(2)}{2}\,({\sqrt 2} \,y(s)z_{2,0}(s) + y(s)^2)
\end{split}
\end{equation*}
By repeating this argument for each group of non-diagonal type the proposition follows.
\end{proof}

We note that these polynomials carry the $L$-spectrum information
only in the variables $x(s)$, $y(s)$, $z_{i,r}(s)$, $1\le i \le
3$. From Propositions \ref{p.zpoly} and \ref{p.zpoly3} the
following result is clear:

\begin{coro}    \label{c.lisosp}
$M_\G$ and $M_{\G'}$ are $L$-isospectral if and only if the
nonzero monomials appearing in  $Z^\G (s)$ and  $Z^{\G'}(s)$ are
the same. In particular, if $\G, \G'$ are of diagonal type  then:
$M_\G$ and $M_{\G'}$ are $L$-isospectral if and only if for each
fixed $1\leq t \le  d\leq 3$ we have
$$\{(t,d) : c_{d,t}(\G)\not= 0\} = \{(t,d) : c_{d,t}(\G')\not= 0\}.$$
\end{coro}

We now give the list of $p$-heat trace polynomials for the groups
in the class considered in this paper, except for  ${\bf 72}, {\bf
74}$, which have holonomy groups of orders 12 and 24 respectively
and can not be $p$-isospectral to any other Bieberbach group.

\pagebreak
\begin{teo}\label{p-heat traces} In the notation above we have:
\sk

\noindent $F\simeq \{\I\}:$
\smallskip

\noindent $ \begin{array}{l} Z_p^{\bf 1}(s)= \tbinom{4}{p}\,
x(s)^4.
\end{array}
$

\noindent $F\simeq \Z_2:$
\smallskip

\noindent
$ \begin{array}{l}
2Z_p^{\bf 2}(s)= \tbinom{4}{p}\, x(s)^4 + K^4_p(1)\; x(s)^2y(s).  \\
2Z_p^{\bf 2'}(s)= \tbinom{4}{p}\, x(s)^4 + K^4_p(1)\; x(s)y(s)^2.  \\
2Z_p^{\bf 2''}(s)= \tbinom{4}{p}\, x(s)^4 + K^4_p(1)\; y(s)^3.  \\
2Z_p^{\bf 3}(s)= \tbinom{4}{p}\, x(s)^4 + \frac {\sqrt 2}2 K^4_p(1)\; x(s)y(s)z_{2,0}(s)= 2Z_p^{\bf 3''}(s).  \\
2Z_p^{\bf 3'}(s)= \tbinom{4}{p}\, x(s)^4 + \frac {\sqrt 2}2 K^4_p(1)\; y(s)^2z_{2,0}(s)= 2Z_p^{\bf 3'''}(s).  \\
%2Z_p^{\bf 3''}(s)= \tbinom{4}{p}\, x(s)^4 + \frac {\sqrt 2}2 K^4_p(1)\; x(s)y(s)z_{2,2}(s).  \\
%2Z_p^{\bf 3'''}(s)= \tbinom{4}{p}\, x(s)^4 + \frac {\sqrt 2}2 K^4_p(1)\; y(s)^2z_{2,2}(s).  \\
2Z_p^{\bf 4}(s)= \tbinom{4}{p}\, x(s)^4 + K^4_p(3)\; y(s).  \\
2Z_p^{\bf 5}(s)= \tbinom{4}{p}\, x(s)^4 + K^4_p(2)\; x(s)y(s).  \\
2Z_p^{\bf 5'}(s)= \tbinom{4}{p}\, x(s)^4 + K^4_p(2)\; y(s)^2.  \\
2Z_p^{\bf 6}(s)= \tbinom{4}{p}\, x(s)^4 + \frac {\sqrt 2}2 K^4_p(2)\; y(s)z_{2,0}(s)=2Z_p^{\bf 6'}(s). \\
%2Z_p^{\bf 6'}(s)= \tbinom{4}{p}\, x(s)^4 + \frac {\sqrt 2}2 K^4_p(2)\; y(s)z_{2,2}(s).
\end{array}  $

\smallskip

\noindent $F\simeq \Z_3:$

\smallskip

\noindent $ \begin{array}{l}
3Z_p^{\bf 47}(s)= \tbinom{4}{p}\,  x(s)^4 + \frac {\sqrt 3}3 tr_p(B)\; (z_{1, 1/3}(s)z_{3,0}(s)+z_{1, 2/3}(s)z_{3,0}(s))=3Z_p^{\bf 47'}(s).  \\
%3Z_p^{\bf 47'}(s)= \tbinom{4}{p}\,  x(s)^4 + \frac {\sqrt 3}3 tr_p(B)\; (z_{1,1}(s)z_{3,1}(s)+z_{1,2}(s)z_{3,2}(s)).
\end{array}  $

\smallskip

\noindent $F\simeq \Z_2^2:$

\smallskip

\noindent $ \begin{array}{l}
4Z_p^{\bf 7}(s)= \tbinom{4}{p}\, x(s)^4 + 2 K^4_p(1)\; x(s)^2y(s) + K^4_p(2)\; x(s)y(s). \\
4Z_p^{\bf 7'}(s)= \tbinom{4}{p}\, x(s)^4 + K^4_p(1)\; (x(s)^2y(s)+x(s)y(s)^2) + K^4_p(2)\; y(s)^2 = 4Z_p^{\bf 10}(s). \\
4Z_p^{\bf 8}(s)= \tbinom{4}{p}\, x(s)^4 + K^4_p(1)\; (x(s)^2y(s) + x(s)y(s)^2) + K^4_p(2)\; x(s)y(s) =4Z_p^{\bf 9'}(s).  \\
4Z_p^{\bf 8'}(s)= \tbinom{4}{p}\, x(s)^4 + K^4_p(1)\; (x(s)^2y(s) + y(s)^3) + K^4_p(2)\; y(s)^2.  \\
4Z_p^{\bf 9}(s)= \tbinom{4}{p}\, x(s)^4 + 2 K^4_p(1)\; x(s)^2y(s) + K^4_p(2)\; y(s)^2.  \\
4Z_p^{\bf 10'}(s)= \tbinom{4}{p}\, x(s)^4 + K^4_p(1)\; (x(s)^2y(s) + y(s)^3) + K^4_p(2)\; x(s)y(s). \\
4Z_p^{\bf 10''}(s)= \tbinom{4}{p}\, x(s)^4 + 2 K^4_p(1)\; x(s)y(s)^2  + K^4_p(2)\; x(s)y(s). \\
4Z_p^{\bf 11}(s)= \tbinom{4}{p}\, x(s)^4 + 2 K^4_p(1)\; x(s)y(s)^2 + K^4_p(2)\; y(s)^2.  \\
4Z_p^{\bf 11'}(s)= \tbinom{4}{p}\, x(s)^4 +  K^4_p(1)\; (x(s)y(s)^2+y(s)^3) + K^4_p(2)\; x(s)y(s).
\end{array}  $

\smallskip

\noindent $ \begin{array}{l}
4Z_p^{\bf 12}(s)= \tbinom{4}{p}\, x(s)^4 + \sqrt 2  K^4_p(1)\; x(s)y(s)z_{2,0}(s) + K^4_p(2)\; y(s)^2. \\
4Z_p^{\bf 12'}(s)= \tbinom{4}{p}\, x(s)^4 + \frac {\sqrt 2}2  K^4_p(1)\; (y(s)^2z_{2,0}(s)+x(s)y(s)z_{2,0}(s)) + K^4_p(2)\; x(s)y(s).
\end{array} $

\smallskip

\noindent $ \begin{array}{l}
4Z_p^{\bf 13}(s)= \tbinom{4}{p}\, x(s)^4 + K^4_p(1)\; (x(s)^2y(s) + \f {\sqrt 2}2 x(s)y(s)z_{2,0}(s)) + \f{\sqrt 2}2  K^4_p(2)\; y(s)z_{2,0}(s). \\
4Z_p^{\bf 13'}(s)= \tbinom{4}{p}\, x(s)^4 + K^4_p(1)\; (y(s)^3 + \f {\sqrt 2}2 x(s)y(s)z_{2,0}(s)) + \f{\sqrt 2}2  K^4_p(2)\; y(s)z_{2,0}(s). \\
4Z_p^{\bf 14}(s)= \tbinom{4}{p}\, x(s)^4 + K^4_p(1)\; x(s)\,(y(s)^2 + \f {\sqrt 2}2 y(s)z_{2,0}(s)) + \f{\sqrt 2}2  K^4_p(2)\; y(s)z_{2,0}(s) = 4Z_p^{\bf 14'}(s). \\
%4Z_p^{\bf 14'}(s)= \tbinom{4}{p}\, x(s)^4 + K^4_p(1)\; (x(s)y(s)^2 + \f {\sqrt 2}2 x(s)y(s)z_{2,2}(s)) + \f{\sqrt 2}2  K^4_p(2)\; y(s)z_{2,0}(s). \\
4Z_p^{\bf 15}(s)= \tbinom{4}{p}\, x(s)^4 + K^4_p(1)\; (x(s)y(s)^2 + \f {\sqrt 2}2 y(s)^2z_{2,0}(s)) + \f{\sqrt 2}2  K^4_p(2)\; y(s)z_{2,0}(s)=4Z_p^{\bf 15'}(s). \\
%4Z_p^{\bf 15'}(s)= \tbinom{4}{p}\, x(s)^4 + K^4_p(1)\; (x(s)y(s)^2 + \f {\sqrt 2}2 y(s)^2z_{2,0}(s)) + \f{\sqrt 2}2  K^4_p(2)\; y(s)z_{2,0}(s). \\
\end{array} $

\smallskip

\noindent $ \begin{array}{l}
4Z_p^{\bf 18}(s)= \tbinom{4}{p}\, x(s)^4 + K^4_p(1)\; x(s)^2y(s) + K^4_p(2)\; x(s)y(s) + K^4_p(3)\; y(s) = 4Z_p^{\bf 19}(s). \\
4Z_p^{\bf 19'}(s)= \tbinom{4}{p}\, x(s)^4 + K^4_p(1)\; x(s)y(s)^2 + K^4_p(2)\; x(s)y(s) + K^4_p(3)\; y(s) = 4Z_p^{\bf 21}(s).\\
4Z_p^{\bf 20}(s)= \tbinom{4}{p}\, x(s)^4 + K^4_p(1)\; x(s)^2y(s) + K^4_p(2)\; y(s)^2 + K^4_p(3)\; y(s). \\
4Z_p^{\bf 20'}(s)= \tbinom{4}{p}\, x(s)^4 + K^4_p(1)\; x(s)y(s)^2 + K^4_p(2)\; y(s)^2 + K^4_p(3)\; y(s). \\
4Z_p^{\bf 21'}(s)= \tbinom{4}{p}\, x(s)^4 + K^4_p(1)\; y(s)^3 + K^4_p(2)\; x(s)y(s) + K^4_p(3)\; y(s). \\
\end{array}  $

\smallskip

\noindent $ \begin{array}{l}
4Z_p^{\bf 22}(s)\!=\! \tbinom{4}{p}\, x(s)^4 + \frac {\sqrt 2}2  K^4_p(1)\; x(s)y(s)z_{2,0}(s)\! +\! \frac {\sqrt 2}2  K^4_p(2)\; y(s)z_{2,0}(s)\! +\!  K^4_p(3)\; y(s)\!=\!4Z_p^{\bf 22'}(s). \\
%4Z_p^{\bf 22'}(s)= \tbinom{4}{p}\, x(s)^4 + \frac {\sqrt 2}2  K^4_p(1)\; x(s)y(s)z_{2,20}(s) + \frac {\sqrt 2}2  K^4_p(2)\; y(s)z_{2,0}(s) +  K^4_p(3)\; y(s).
\end{array} $

\smallskip

\noindent $ \begin{array}{l}
4Z_p^{\bf 23}(s)= \tbinom{4}{p}\, x(s)^4 + K^4_p(2)\; x(s)y(s) + 2 K^4_p(3)\; y(s). \\
4Z_p^{\bf 23'}(s)= \tbinom{4}{p}\, x(s)^4 + K^4_p(2)\; y(s)^2 + 2 K^4_p(3)\; y(s).
\end{array} $

\smallskip

\noindent $ \begin{array}{l}
4Z_p^{\bf 24}(s)= \tbinom{4}{p}\, x(s)^4 + 3 K^4_p(2)\; x(s)y(s) = 4Z_p^{\bf 26}(s).\\
4Z_p^{\bf 25}(s)= \tbinom{4}{p}\, x(s)^4 + K^4_p(2)\; (2 x(s)y(s) + y(s)^2). \\
4Z_p^{\bf 27}(s)= \tbinom{4}{p}\, x(s)^4 + K^4_p(2)\; (x(s)y(s) + 2 y(s)^2).  \\
\end{array}  $

\smallskip

\noindent $ \begin{array}{l}
4Z_p^{\bf 28}(s)= \tbinom{4}{p}\, x(s)^4 +  K^4_p(2)\; (x(s)y(s) + \sqrt 2 \, y(s)z_{2,0}(s)). \\
\end{array} $

\smallskip

\noindent $ \begin{array}{l}
4Z_p^{\bf 29}(s)= \tbinom{4}{p}\, x(s)^4 +  K^4_p(2)\; (y(s)^2 + \sqrt 2 \, y(s)z_{2,0}(s)). \\
4Z_p^{\bf 29'}(s)= \tbinom{4}{p}\, x(s)^4 +  K^4_p(2)\; (x(s)y(s)
+ \frac {\sqrt 2}2  y(s)z_{2,0}(s) + \frac {\sqrt 2}2
y(s)z_{2,1/2}(s)).
\end{array} $

\smallskip

\noindent $F\simeq \Z_4:$

\smallskip

\noindent $ \begin{array}{l}
4Z_p^{\bf 45}(s)= \tbinom{4}{p}\, x(s)^4 +  K^4_p(2)\; x(s)y(s) + tr_p(B)\; x(s)(z_{1,1/4}(s) + z_{1,3/4}(s)). \\
4Z_p^{\bf 45'}(s)= \tbinom{4}{p}\, x(s)^4 +  K^4_p(2)\; y(s)^2 + tr_p(B)\; (z_{1,1/4}(s)^2 + z_{1,3/4}(s)^2). \\
4Z_p^{\bf 50}(s)= \tbinom{4}{p}\, x(s)^4 +  K^4_p(2)\; x(s)y(s) + tr_p(B)\; (z_{1,1/4}(s) + z_{1,3/4}(s)). \\
4Z_p^{\bf 51}(s)= \tbinom{4}{p}\, x(s)^4 +  K^4_p(2)\; y(s)^2 +
\sqrt 2 \, tr_p(B)\; z_{2,1/2}(s).
\end{array} $

\smallskip

\noindent
 $F\simeq \Z_6:$

\smallskip

\noindent $ \begin{array}{cl}
6Z_p^{\bf 64}(s)  = & \!\!\!\tbinom{4}{p}\, x(s)^4 + \frac {\sqrt 3}3  tr_p(B^2)\;
(z_{1,1/3}(s)z_{3,0}(s)+z_{1,2/3}(s)z_{3,0}(s)) +  K^4_p(3)\; y(s) \,+ \\
  &   tr_p(B)\; (z_{1,1/6}(s) + z_{1,5/6}(s)).
\end{array} $

\smallskip

\noindent $F\simeq D_3:$

\smallskip
\noindent $ \begin{array}{cl} 6Z_p^{\bf 67}(s)= &\!\!\! \tbinom{4}{p}\,
x(s)^4 + 3\frac{\sqrt 2}2 K^4_p(2)\, y(s)z_{2,0}(s)+
  \frac {\sqrt 3}3  tr_p(B)\;
(z_{1,1/3}(s)+z_{1,2/3}(s)) z_{3,0}(s).
\end{array} $

\smallskip

%\pagebreak

\noindent $F\simeq \Z_2^3:$

\smallskip

\noindent $ \begin{array}{l} 8Z_p^{\bf 33}(s)=\tbinom{4}{p}\,
x(s)^4 + 3 K^4_p(1)\; x(s)^2y(s)
+ K^4_p(2)\; (2x(s)y(s)+y(s)^2) +  K^4_p(3)\; y(s). \\
8Z_p^{\bf 34}(s)=\tbinom{4}{p}\, x(s)^4\!+\!K^4_p(1)\, (2x(s)^2y(s) + x(s)y(s)^2)\!+\!K^4_p(2)\, (x(s)y(s)+2y(s)^2)\!+\!K^4_p(3)\; y(s). \\
8Z_p^{\bf 35}(s)=\tbinom{4}{p}\, x(s)^4 + K^4_p(1)\; (x(s)^2y(s)+x(s)y(s)^2+y(s)^3) +  3 K^4_p(2)\; x(s)y(s) + K^4_p(3)\; y(s).  \\
8Z_p^{\bf 36}(s)= \tbinom{4}{p}\, x(s)^4 + K^4_p(1)\; (2x(s)^2y(s) + x(s)y(s)^2) + 3 K^4_p(2)\; x(s)y(s) +  K^4_p(3)\; y(s).\\
8Z_p^{\bf 37}(s)= \tbinom{4}{p}\, x(s)^4 +  K^4_p(1)\; (2x(s)^2y(s)+y(s)^3) + K^4_p(2)\; (2x(s)y(s)+y(s)^2) +  K^4_p(3)\; y(s).\\
8Z_p^{\bf 38}(s)= 8Z_p^{\bf 34}(s). \\
%\tbinom{4}{p}\, x(s)^4 + K^4_p(1)\; (2x(s)^2y(s) + x(s)y(s)^2) + K^4_p(2)\; (x(s)y(s)+2y(s)^2) + K^4_p(3)\; y(s)\\
8Z_p^{\bf 39}(s)= \tbinom{4}{p}\, x(s)^4\! +\!  K^4_p(1)\, (x(s)^2y(s)+2x(s)y(s)^2)\! +\! K^4_p(2)\, (2x(s)y(s)+y(s)^2)\! +\!  K^4_p(3)\, y(s). \\
8Z_p^{\bf 40}(s)= 8Z_p^{\bf 35}(s).\\
8Z_p^{\bf 41}(s)= 8Z_p^{\bf 39}(s). \\
8Z_p^{\bf 42}(s)= \tbinom{4}{p}\, x(s)^4 +  3 K^4_p(1)\; x(s)y(s)^2 + 3 K^4_p(2)\; x(s)y(s) +  K^4_p(3)\; y(s).\\
\end{array} $

\smallskip

\noindent $ \begin{array}{l}
8Z_p^{\bf 43}(s)= \tbinom{4}{p}\, x(s)^4 +  K^4_p(1)\; x(s)^2y(s) +  K^4_p(2)\; (2x(s)y(s)+y(s)^2) + 3 K^4_p(3)\; y(s).\\
8Z_p^{\bf 44}(s)= \tbinom{4}{p}\, x(s)^4 +  K^4_p(1)\; x(s)y(s)^2 +  3 K^4_p(2)\; x(s)y(s) + 3 K^4_p(3)\; y(s).\\
\end{array}        $

\smallskip

\noindent
 $F\simeq \Z_2 \times \Z_4:$

\smallskip

\noindent $ \begin{array}{cl} 8Z_p^{\bf 57}(s) = & \!\! \!\tbinom{4}{p}\,
x(s)^4 + K^4_p(1)\; x(s)y(s)^2 + K^4_p(2)\; x(s)y(s) +
tr_p(B_1)\; (z_{1,1/4}(s)+z_{1,3/4}(s))\, + \\
 & \,tr_p(B_1B_2)\; x(s)(z_{1,1/4}(s) + z_{1,3/4}(s)) + K^4_p(3)\; y(s)
  = 8Z_p^{\bf 58}(s).
%\\8Z_p^{\bf 58}(s) =\tbinom{4}{p}\, x(s)^4 +  K^4_p(1)\; x(s)y(s)^2 +   K^4_p(2)\; x(s)y(s) + \\  tr_p(B_1B_2)\; (y(s)z_{1,2}(s)+y(s)z_{1,6}(s))    +\;\; tr_p(B_1)\; (z_{1,2}(s)+z_{1,6}(s)) +  K^4_p(3)\; y(s).\\
\end{array} $

\smallskip

\noindent $F\simeq D_4:$

\smallskip

\noindent $ \begin{array}{cl} 8Z_p^{\bf 54}(s) =& \!\!\! \tbinom{4}{p}\,
x(s)^4 + \sqrt 2 \, K^4_p(1)\; x(s)y(s)z_{2,0}(s) +
  K^4_p(2)\; (y(s)^2 + \, {\sqrt 2}\,  y(s)z_{2,0}(s))\; + \\&
\sqrt 2\, tr_p(B_1)\; z_{2,1/2}(s).
\end{array} $

\smallskip

\noindent $ \begin{array}{cl} 8Z_p^{\bf 56}(s) =& \!\!\! \tbinom{4}{p}\,
x(s)^4 + K^4_p(2)\; y(s)(x(s) + {\sqrt 2}\, z_{2,0}(s))\, + tr_p(B_1)\;
(z_{1,1/4}(s)+z_{1,3/4}(s))\; + \\ & 2 K^4_p(3)\; y(s).
\end{array} $
\end{teo}

\smallskip

\noindent $ \begin{array}{cl} 8Z_p^{\bf 60}(s)  =&\!\!\!
\tbinom{4}{p}\, x(s)^4  \!+\!   K^4_p(2)\; (3 x(s)y(s) \!+\!
\sqrt 2 \, y(s)z_{2,0}(s))\!+\! tr_p(B_1)\; x(s) (z_{1,1/4}(s) + z_{1,3/4}(s)).  \\
8Z_p^{\bf 61}(s)  =& \!\!\! \tbinom{4}{p}\, x(s)^4   + K^4_p(2)\;
 (x(s)y(s) + 2y(s)^2 +   {\sqrt 2}\, y(s)z_{2,0}(s))\;  +
\\ &  tr_p(B_1)\; x(s) (z_{1,1/4}(s) +  z_{1,3/4}(s)).  \\
8Z_p^{\bf 62}(s) =&\!\!\! \tbinom{4}{p}\, x(s)^4 + K^4_p(2)\; (3 x(s)y(s)
+ {\sqrt 2}\, y(s)z_{2,0}(s))\, +\! tr_p(B_1)\; y(s)(z_{1,1/4}(s) +
z_{1,3/4}(s)).
\end{array} $

\bigskip
\begin{rem}
(a) With a similar asymptotic argument as the one used in Lemma
\ref{algindep}, one can show that the functions $x(s), y(s),
z_{i,k}(s)$ for $1\le i \le 3,\, 0\le k \le |F|-1$, appearing in
the previous theorem are algebraically independent. This fact will
be used in the proof of the main results in the following section.

(b) Let $\G$ be a 4-dimensional Bieberbach group with canonical
lattice and holonomy group $F$. We note that one can read off some
properties of $\Gamma$ from the expression of the zeta function.
For instance: {\em (i)} $\G$ has diagonal holonomy representation
if and only if $Z_p^\G (s)\in \tfrac 1{|F|} \Z[x(s),y(s)]$.
\noindent {\em (ii)} $M_\G$ is orientable if and only if for every
$i,j,k$ with $i+j+k$ odd, the coefficient of $x(s)^iy(s)^jz(s)^k$
in $Z_0^\G (s)$ is zero.
\end{rem}

\section{Isospectral classifications}           \label{s.sec4}
    In this section we shall use the information on the
    $p$-heat trace polynomials in the previous sections to
determine all isospectral classes  for the different type of
spectra for all compact flat manifolds of dimension 4 in the class
considered. We will thus find all pairs that are isospectral,
$p$-isospectral, $L$-isospectral or $[L]$-isospectral in this
case.

\begin{teo}[$p$-isospectrality]     \label{t.pisosp}
Let $\G, \G'$ be two 4-dimensional Bieberbach groups with canonical translation lattice. Then:
\begin{enumerate}
    \item [{\em (i)}] %\underline
{If $p=0,4$}: $M_\G$ and ${M_{\G'}}$ are $p$-isospectral if and only if
$\{\G, \G'\}$ is one of the following: $\{{\bf 3},{\bf 3''}\}$, $\{{\bf
3'},{\bf 3'''}\}$, $\{{\bf 6},{\bf 6'}\}$ ($F\simeq \Z_2$), $\{{\bf
47},{\bf 47'}\}$ ($F\simeq \Z_3$), $\{{\bf 7'},{\bf 10}\}$, $\{{\bf
8},{\bf 9'}\}$, $\{{\bf 14},{\bf 14'}\}$, $\{{\bf 15},{\bf 15'}\}$,
$\{{\bf 18},{\bf 19}\}$, $\{{\bf 19'},{\bf 21}\}$, $\{{\bf 22},{\bf
22'}\}$, $\{{\bf 24},{\bf 26}\}$ ($F\simeq \Z_2^2$), or $\{{\bf 34},{\bf
38}\}$, $\{{\bf 35},{\bf40}\}$, $\{{\bf 39},{\bf41}\}$ ($F\simeq \Z_2^3$),
$\{{\bf 57},{\bf 58}\}$ ($F\simeq \Z_2\times \Z_4$).
    \item [{\em (ii)}]
{If $p=1,3$}: $M_\G$  and ${M_{\G'}}$ are $p$-isospectral if and only if $\G, \G'$ both belong to one of the following families: %\linebreak
 $\{{\bf 3},{\bf 3''}\}$,
$\{{\bf 3'},{\bf 3'''}\}$, $\{{\bf 5},{\bf 5'},{\bf 6},{\bf 6'}\}$
($F\simeq \Z_2$), $\{{\bf 47},{\bf 47'}\}$ ($F\simeq \Z_3$), $\{{\bf
7},{\bf 9}\}$, $\{{\bf 7'},{\bf 8},{\bf 9'},{\bf 10}\}$, $\{{\bf 8'}, {\bf
10'}\}$, $\{{\bf 10''},{\bf 11}\}$, $\{{\bf 14},{\bf 14'}\}$, $\{{\bf
15},{\bf 15'}\}$, $\{{\bf 18},{\bf 19},{\bf 20}\}$, $\{{\bf 19'},{\bf
20'},{\bf 21}\}$, $\{{\bf 22},{\bf 22'}\}$, $\{{\bf 23},{\bf 23'}\}$
($F\simeq \Z_2^2$), $\{{\bf 24},{\bf 25},{\bf 26},{\bf 27},{\bf 28},{\bf
29},{\bf 29'},{\bf 50},{\bf 51}\}$ ($F \simeq \Z_2^2, \Z_4$), $\{{\bf
34}$, ${\bf 36}, {\bf 38}\}$, $\{{\bf 35},{\bf 40}\}$, $\{{\bf 39},{\bf
41}\}$ ($F\simeq \Z_2^3$), $\{{\bf 57},{\bf 58}\}$ ($F\simeq \Z_2\times
\Z_4$) or $\{{\bf 60},{\bf 61}\}$ ($F\simeq D_4$).
    \item [{\em (iii)}]
{If $p=2$}: $M_\G$  and ${M_{\G'}}$ are $2$-isospectral if and
only if $\G, \G'$ both belong to one of the following families:
 $\{{\bf 2,2',2'',3,3',3'',3''',4}\}$, $\{{\bf 6},{\bf 6'}\}$
 ($F\simeq\Z_2$), $\{{\bf 47},{\bf 47'}\}$ ($F\simeq \Z_3$),
 $\{{\bf 7,8,9',10',10'',11',12',18,19,19',21,21',23,50}\}$,
 $\{{\bf 7',8',9,10,11,12,20,20',23'}$, \linebreak ${\bf 51}\}$ ($F \simeq
  \Z_2^2,\Z_4$),  $\{{\bf 13,\bf 13',\bf 14, \bf 14', \bf 15, \bf 15'},
  {\bf 22,22'}\}$, $\{{\bf 24,26}\}$ ($F\simeq\Z_2^2$),
   $\{{\bf 33,37,39}$, ${\bf 41},{\bf 43}\}$, $\{{\bf 34,38}\}$,
   $\{{\bf 35,36,40,42,44}\}$
   ($F\simeq \Z_2^3$), $\{{\bf 57},{\bf 58}\}$
($F\simeq \Z_2\times \Z_4$).
\end{enumerate}
\end{teo}

\begin{proof}
By Theorem \ref{t.zetafunc}, if $M_\G$ and ${M_\G}'$ are
$p$-isospectral then $|F|=|F'|$. Therefore when looking for
$p$-isospectrality  we only need to compare pairs both  with
$F\simeq \Z_2, F\simeq \Z_3$, $F\simeq \Z_2^2$ or $\Z_4$, $F\simeq
\Z_2^3, \Z_2 \times \Z_4$ or $D_4$ or else $F\simeq \Z_6$ or
$D_3$. Also, for each $p$, $\beta_p$ gives the multiplicity of the
eigenvalue 0 of $\Delta_p$, hence we need only consider pairs
having the same $\beta_p$.

Take $|F|=2$. Here the groups are ${\bf
2,2',2'',3,3',3'',3''',4,5,5'}$, ${\bf 6, 6'}$. Looking at the
corresponding $p$-heat trace polynomials in Section 3, we see that
they are all different from each other, hence they can not be
pairwise isospectral. However, the vanishing of the Krawtchouk
polynomials for some values of $p$ yields some $p$-isospectrality
for $1\leq p\leq 3$. For instance, for $p=1$ or 3 we know that
$K_p^4(2)=0$, so $Z_p^{\bf 5}(s)=  Z_p^{\bf 5'}(s)= Z_p^{\bf
6}(s)=  Z_p^{\bf 6'}(s)= 2\, x(s)^4$. Thus ${\bf 5}$, ${\bf 5'}$,
${\bf 6}$ and ${\bf 6'}$ are 1-isospectral and 3-isospectral. In
the same way, $K_p^4(1)=K_p^4(3)=0$ for $p=2$ and so $Z_2(s)=3\,
x(s)^4$ for all groups ${\bf 2, 2', 2'', 3, 3', 3'', 3''', 4}$,
hence these groups %${\bf 2, 2', 2'', 3, 3', 3'', 3''', 4}$
are 2-isospectral to each other.

In the cases of the remaining groups, i.e.\@ those with $|F|=3,
4,6$ and $8$, if $p=0$, using the complete list of the $p$-heat
trace polynomials  given in Section~3, we easily see that only the
pairs $\{{\bf 3},{\bf 3''}\}$, $\{{\bf 3'},{\bf 3'''}\}$, $\{{\bf
6},{\bf 6'}\}$, $\{{\bf 47},{\bf 47'}\}$, $\{{\bf 7'},{\bf 10}\}$,
$\{{\bf 8},{\bf 9'}\}$, $\{{\bf 14},{\bf 14'}\}$, $\{{\bf 15},{\bf
15'}\}$, $\{{\bf 18},{\bf 19}\}$, $\{{\bf 19'},{\bf21}\}$, $\{{\bf
22},{\bf 22'}\}$, $\{{\bf 24},{\bf26}\}$,  $\{{\bf 34},{\bf38}\}$,
$\{{\bf 35},{\bf40}\}$, $\{{\bf 39},{\bf41}\}$, $\{{\bf 57},{\bf
58}\}$ have the same $p$-heat trace expressions for $p=0$.
Actually they have the same $p$-heat trace for all $p$, thus they
are $p$-isospectral for $0\le p \le 4$.

In the case of the $p$-spectrum for $1\leq p\leq 3$,  again by comparison of the
$p$-heat traces one obtains the $p$-isospectral sets asserted in the theorem.
\end{proof}

\begin{coro} Let $\G,\G'$ be as in Theorem~\ref{t.pisosp}. Then
$M_\G$ and ${M_\G}'$ are isospectral if and only if they are $p$-isospectral for
 $0 \le p \le 4$.
\end{coro}

\begin{rem}
(i) Observe that the groups in each isospectral pair in
Theorem~\ref{t.pisosp} have the same holonomy representation.

(ii) The assertion in the corollary fails to be true for general
flat manifolds. Examples of $0$-isospectral and not
$1$-isospectral flat manifolds in dimension $n\ge 6$ are given in
\cite{MR2}. In \cite{RC} two nonhomeomorphic flat 3-manifolds that
are isospectral on functions are given, showing this is the only
such pair. These manifolds have different holonomy groups: $\Z_4$
and $\Z_2\times \Z_2$ respectively and one can check they are not
isospectral on $1$-forms by comparison of the $p$-heat traces.

(iii) The theorem also shows that isospectrality on functions is
much less common than  $p$-isospec\-tra\-lity for $1\leq p \leq
3$. Some  facts can be observed. For $p=1,3$, one of the
$p$-isospectral pairs, $\{{\bf 60,61}\}$, involves two groups with
non-abelian holonomy group $D_4$. Also, in the $p$-isospectral set
$\{{\bf 24,25,26,27,28,29,29',50,51}\}$  five different holonomy
representations are present. Here,  the  first seven groups have
holonomy  $F\simeq \Z_2^2$ and give orientable manifolds, while
the  last two groups have holonomy $F\simeq \Z_4$ and produce
non-orientable manifolds.  Relative to $2$-isospectral classes,
again some manifolds are orientable and others not, and different
holonomy groups occur.
\end{rem}
\medskip

We now turn into the study of length spectra for the flat
$4$-manifolds considered in this paper. It is a classical result
that in the case of flat tori the eigenvalue spectrum and the
length spectrum determine each other. This is not true for general
Riemannian manifolds, even in the flat case. There exist
$4$-dimensional flat  manifolds (namely ${\bf 8},{\bf 9'}$ in our
notation) which are isospectral but such that some multiplicities
of closed geodesics are different (see \cite{MR4}, Ex.\@ 3.4, up
to isometry; see also the references therein for earlier examples
in the case of nilmanifolds). In particular, this shows that in
general, the eigenvalue spectrum does not determine the length
spectrum. In the converse direction, in \cite{MR4} an example was
given, for $n=13$, of $[L]$-isospectral flat manifolds which are
not isospectral on functions.

The study of $[L]$-isospectrality is delicate since multiplicities must be taken into account. We recall that the  multiplicity of a length $l$ is the number of conjugacy classes in $\G$ having length $l$. We shall see that already in the context in this paper, most $L$-isospectral pairs  fail to be $[L]$-isospectral.
The next result describes the situation.

\begin{teo}[$L$-isospectrality]  \label{t.Lisosp}Let $\G, \G'$ be two 4-dimensional
Bieberbach groups with canonical translation lattice.

(i) The $L$-isospectral sets are the following: $\{{\bf 3},{\bf
3''}\}$, $\{{\bf 3'},{\bf 3'''}\}$, $\{{\bf 6},{\bf 6'}\}$,
$\{{\bf 47},{\bf 47'}\}$, $\{{\bf 7'},{\bf 10}\}$, $\{{\bf 8},{\bf
9'}\}$, $\{{\bf 14},{\bf 14'}\}$, $\{{\bf 15},{\bf
15'}\}$,\,$\{{\bf 18},{\bf 19}\}$,\,$\{{\bf
19'},{\bf21}\}$,\,$\{{\bf 22},{\bf 22'}\}$, $\{{\bf
24},{\bf26}\}$,\,$\{{\bf 34},{\bf 38},{\bf 39},
{\bf41}\}$,\,$\{{\bf 35},{\bf40}\}$, $\{{\bf 57},{\bf 58}\}$ and,
in addition, the sets  $\{{\bf 25},{\bf 27}\}$, $\{{\bf 33},{\bf
43}\}$, $\{{\bf 42},{\bf 44}\}$.

(ii)  The $[L]$-isospectral pairs are $\{{\bf 3},{\bf 3''}\}$,
$\{{\bf 3'},{\bf 3'''}\}$, $\{{\bf 6},{\bf 6'}\}$, $\{{\bf
47},{\bf 47'}\}$, $\{{\bf 7'},{\bf 10}\}$, $\{{\bf 14},{\bf
14'}\}$, $\{{\bf 15},{\bf 15'}\}$, $\{{\bf 18},{\bf 19}\}$,
$\{{\bf 19'},{\bf21}\}$, $\{{\bf 22}, {\bf 22'}\}$, $\{{\bf
24},{\bf26}\}$,  and $\{{\bf 57},{\bf 58}\}$.
\end{teo}

\begin{proof}
(i) All isospectral pairs in (i) of Theorem~\ref{t.pisosp} are
$L$-isospectral. As we have seen in Section~3, the information on
the $L$-spectrum of $\G$ depends only on the monomials ocurring in
the $0$-heat trace. Since for the $L$-spectrum we do not look at
multiplicities, $\G,\G'$ will be $L$-isospectral if and only if
the same monomials occur in the $0$-heat trace polynomials (see
corollary~\ref{c.lisosp}).

Checking the list of $p$-heat traces  we see, for instance, that
$4Z_p^{\bf 25}(s)=\tbinom{4}{p}\,x(s)^4 + K_p^4(2)\,\big(y(s)^2 +
2\,x(s)y(s)\big)$ and $4Z_p^{\bf 27}(s)=\tbinom{4}{p}\,x(s)^4 +
K_p^4(2)\,\big(x(s)y(s) + 2\,y(s)^2\big)$, thus ${\bf 25}$ and
${\bf 27}$ are $L$-isospectral.  In the same way we see that
$\{{\bf 33},{\bf 43}\}$, $\{{\bf 42},{\bf 44}\}$ and $\{{\bf
34,38, 39,41}\}$ are $L$-isospectral.

\smallskip

(ii) Since $[L]$-isospectral implies $L$-isospectral, we must only
consider the sets given in (i).

In \cite{MR4}, the length spectrum of flat manifolds is studied
and a criterion to prove $[L]$-isospectrality is given (see
Proposition 3.1). This criterion depends on the existence of a
bijection $\Phi$ with certain properties  between partitions
${\mathcal P}$ and ${\mathcal P}'$ of $F$ and $F'$ respectively.
Also,  in \cite{MR4} the criterion is applied to show
$[L]$-isospectrality of some pairs, including the pair $\{{\bf 24,
26}\}$ (see Example 3.3, up to isomorphism).

We shall sketch the verification for the pair $\{{\bf 3},{\bf
3''}\}$ (see Table 2.1). We have that $\G$ (resp.\@ $\G''$) is
generated by $\g = d(I,J)L_b$ (resp.\@ $\g'' = d(I,J)L_{b''}$) and
$L_\Ld$, where $b={\frac{e_1}2}$, (resp.\@ $b'' = \frac{e_1+e_3
+e_4}2$).

The elements of $\G$ (resp.\@ $\G''$) are either of the form $\g
L_\ld$ (resp.\@ $\g'' L_\ld$) or $L_\ld$ with $\ld = m_1 e_1 +
m_2e_2 + m_3 e_3 + m_4 e_4 \in \Ld$, i.e.\@ $m_i \in \Z$.

We define a correspondence $\phi$ from $\G$ to $\G''$ such that
$L_\ld \mapsto L_\ld$ and $\g L_\ld \mapsto \g'' L_\ld L_{-e_3}$.
We note that both of these elements have squared length $\big(m_1+
\f 12\big)^2 + m_2^2 + \tf12 (m_3+ m_4)^2$.

If we use the relation (\ref{conj}) below we see that $\g L_\ld
\sim \g L_\ld L_{m(e_3 - e_4)}$  and similarly $\g'' L_\ld \sim
\g'' L_\ld L_{m(e_3 - e_4)}$, for any $m \in \Z$. We note that
(\ref{conj2}) does not introduce new conjugacy relations among the
elements of $\G$ and $\G''$, since the holonomy group is cyclic.

It is now easy to see that the map $\phi$ induces the bijection
$\Phi= \I$ from $F$ to $F''=F$, satisfying all the conditions in
Proposition 3.1 of \cite{MR4}, hence $\{{\bf 3},{\bf 3''}\}$ are
$[L]$-isospectral.

The $[L]$-isospectrality of the remaining pairs in (ii) of the
theorem, can be proved similarly, with $\Phi= \I$ as the
bijection. We shall omit this verification.

\msk To conclude the proof of the theorem we must show that the
remaining $L$-isospectral sets, $\{{\bf 8},{\bf 9'}\}$, $\{{\bf
25}, {\bf 27}\}$, $\{{\bf 33},{\bf43}\}$, $\{{\bf 34,38,
39,41}\}$, $\{{\bf 35},{\bf40}\}$, $\{{\bf 42},{\bf 44}\}$, are
not $[L]$-isospectral.

The pair $\{{\bf 8},{\bf 9'}\}$ was shown not to be
$[L]$-isospectral in \cite{MR4} (up to isometry, see Example 3.4).
We now discuss in some detail the case $\{{\bf 25,27}\}$.

If $l$ is the length of a closed geodesic in $M_\G$, the multiplicity is given by
$m_\G(l)=\sum_{\g\in \Ld\backslash \G} m_\g(l)$  where
$$ m_\g(l)=\#\{[\g L_\ld]\in [\G]:\ld \in \Ld,\; \ell(\g L_\ld)=l\}.$$
In order to compute these multiplicities it is necessary to have a
parametrization of the $\G$-conjugacy classes in $\G$. This is
complicated in general but it simplifies for  Bieberbach groups of
diagonal type. From \cite{MR4} we have the following relations.
For general $\G$, if $\g_i=B_iL_{b_i}, \g_j=B_jL_{b_j} \in \G$ and
$\ld, \mu \in \Ld$,the conjugations $\g_iL_\ld{\g_i}^{-1}$ and
$L_\mu(\g_iL_\ld)L_{-\mu}$ give the relations:
\begin{equation} \label{conj}
L_\ld \sim L_{B_i\ld} \qquad \text{ and } \qquad
\g_iL_\ld \sim \g_iL_\ld L_{({B_i}^{-1}-\I)\mu}
\end{equation}
Furthermore, if $F$ is abelian,
the conjugation $\g_j (\g_iL_\ld) {\g_j}^{-1}$, with $i\not =j$,  gives the relation
\begin{equation}  \label{conj2}
\g_iL_\ld \sim \g_i L_{\eta_{j,i}(\ld)}
 \qquad \text{ where } \qquad \eta_{j,i}(\ld):= B_j\ld + (B_j-\I)b_i + B_j(B_i-\I)b_j
\end{equation}

The next table gives the elements $\g_1, = B_1b_1, \g_2 = B_2
Lb_2, \g_{12} = B_{12}b_{12}$ in {\bf 25, 27}, in column notation,
i.e. showing  in the columns, for each $BL_b\in \Gamma$, the
diagonal entries of $B$ together with the coordinates of the
corresponding translation vector $b$ (see \cite{MR2}).

\renewcommand{\arraystretch}{1}
\medskip
\begin{center}
 \begin{tabular}{|rll|rll|rll|}  \hline   $B_1$  &
 $b_1^{25}$ & $b_1^{27}$& $B_2$  & $b_2^{25}$ & $b_2^{27}$ & $B_1B_2$ &$b_{12}^{25}$
 & $b_{12}^{27}$  \\ \hline
-1 && & 1 &{\scriptsize 1/2} &{\scriptsize 1/2} & -1 &{\scriptsize 1/2}&{\scriptsize 1/2} \\
-1 && & -1 &{\scriptsize 1/2} &{\scriptsize 1/2} & 1 &{\scriptsize 1/2} &{\scriptsize 1/2}  \\
1 & & {\scriptsize 1/2} &  -1 && & -1 & & {\scriptsize 1/2} \\
1 & {\scriptsize 1/2} & & 1 & &{\scriptsize 1/2} & 1 &{\scriptsize 1/2} &{\scriptsize 1/2}  \\
 \hline
\end{tabular}
\end{center}

The squared lengths of closed geodesics corresponding to elements not in $\Ld$ are given by
\begin{center} \begin{tabular}{lccc}
${\bf 25}:$ & $m_3^2 + (m_4+\frac12)^2$, & $(m_1+\frac12)^2 + m_4^2$, &  $(m_2+\frac12)^2 + (m_4+\frac12)^2$,  \\
${\bf 27}:$ & $(m_3+\frac12)^2 + m_4^2$, & $(m_1+\frac12)^2 + (m_4+\frac12)^2$, &  $(m_2+\frac12)^2 + (m_4+\frac12)^2$,  \\
\end{tabular} \end{center}
with $m_i\in \Z$.

\sk The  minimal length in both cases is $l=\frac12$. We shall
show that this length has different  multiplicies for ${\bf25}$
and ${\bf27}$. In $\bf 25$ the elements with $l=\frac12$ have the
form:
$$ B_1L_{\frac{e_4}{2}+\ld_1} \text{ where }  \ld_1=m_1e_1+m_2e_2-m_4e_4 \quad \text{ with } m_1,m_2 \in \Z, \;m_4\in \{0,1\}$$
$$ B_2L_{\frac{e_1}{2}+\ld_2} \text{ where }  \ld_2=-n_1e_1+n_2e_2+n_3e_3 \quad \text{ with }
  n_2,n_3 \in \Z, \;n_1\in \{0,1\}.$$
In ${\bf 27}$ the  elements with $l=\frac12$ have the form:
$$ B_1L_{\frac{e_3}{2}+{\ld_1}'} \text{ where }  {\ld_1}'=r_1e_1+r_2e_2-r_3e_3 \quad \text{ with }r_1,r_2 \in \Z, \;r_3\in \{0,1\}.$$

It follows from (\ref{conj}) that $m_1, m_2, n_2,n_3, r_1, r_2$ can be
taken$\mod(2)$.

Using the  relations (\ref{conj2}) we have:
 \begin{align*}
\g_1 L_{m_1e_1+m_2e_2+m_4e_4} & \sim
\g_1 L_{(m_1-1)e_1-(m_2+1)e_2+m_4e_4}  \qquad & \text{ in }{\bf25}& \sk\\
\g_2 L_{n_1e_1+n_2e_2+n_3e_3} & \sim \g_2
L_{-(n_1+1)e_1-(n_2+1)e_2+n_3e_3} \qquad & \text{ in }{\bf 25}&
\sk \\
\g_1L_{r_1e_1+r_2e_2+r_3e_3}  & \sim
\g_1L_{(r_1-1)e_1-(r_2+1)e_2-r_3e_3} \qquad  & \text{ in }{\bf
27}&.
\end{align*}

These three relations divide by 2 the number of relations, and
there are no other relations. Thus $m_{{\bf25}}(\frac12)=4+4=8$
while $m_{{\bf27}}(\frac12)=4$. This shows that ${\bf25},{\bf27}$
do not have the same multiplicities.

The verification of non $[L]$-isospectrality in the remaining
cases can be done similarly as for {\bf 25}, {\bf 27} (by
comparing the multiplicities of small lengths) and will be
omitted.
\end{proof}
\noindent

%                   ----------------
%                   THE BIBLIOGRAPHY
%                   ----------------

%                   End THE BIBLIOGRAPHY

\end{document}